\documentclass[final]{siamltex}

\usepackage{amsmath}
\usepackage{amssymb}
\usepackage{amstext}

\usepackage{graphicx}
\usepackage{hhline}
\usepackage{psfrag}
\usepackage{float}
\usepackage[dvips]{epsfig}
\usepackage{subfig}

\newcommand{\R}{\mathrm{I\hspace{-0.5ex}R}}

\newcommand{\eps}{\epsilon}

\newtheorem{remark}{Remark}

\newtheorem{algorithm}{Algorithm}
\newtheorem{example}{Example}

\begin{document}

\title{A multi-scale  particle   method for mean field equations: the general case}

\author{A. Klar \footnotemark[1]
  \footnotemark[2] \and  S. Tiwari  \footnotemark[1] }
\footnotetext[1]{Technische Universit\"at Kaiserslautern, Department of Mathematics, Erwin-Schr\"odinger-Stra{\ss}e, 67663 Kaiserslautern, Germany 
  (\{klar, tiwari\}@mathematik.uni-kl.de)}
\footnotetext[2]{Fraunhofer ITWM, Fraunhoferplatz 1, 67663 Kaiserslautern, Germany}


\maketitle

\begin{abstract}
A multi-scale meshfree particle method  for macroscopic mean field approximations  of generalized interacting particle models is developed and investigated. The method is working in a uniform way  for   large and small interaction radii.
The well resolved case for large interaction radius is treated, as well as  underresolved situations with  small values of the interaction radius.
In the present work we extend the approach from \cite{KT15} for porous media type limit equations to a more general case, including in particular hyperbolic limits.
The method  can be viewed as  a numerical  transition between a DEM-type method for microscopic  interacting particle systems  and a  meshfree particle method for macroscopic equations.
We discuss in detail the numerical performance of the scheme for various examples and the potential gain in computation time.
The latter is shown to be particularly high for situations near the macroscopic limit.
There are various applications of the method to problems involving mean field approximations in swarming, traffic, pedestrian or granular flow simulation.
\end{abstract}




\section{Introduction}

Interacting particle models are used in many applications ranging from engineering sciences to social sciences and biology. We refer to
\cite{CFRT,DM1,Deg13,Hel95,HS09,HS11} for recent work. 
In certain limits, for   large  numbers 
  of particles, the solution of the microscopic model can be approximated by the solution of so-called mean field or macroscopic equations.
  For a mathematical investigation of these  limits we refer to
\cite{Golse,BH,CCR}. 
Microscopic interacting particle models require numerically the solution of large systems of ordinary differential equations. Many different numerical approaches can be employed for the numerical solution of the macroscopic approximations. In particular, meshless or particle methods are a popular  way to solve these problems.
We refer to \cite{BGLX00,MR87} for the analysis of particle methods, to  \cite{GM77} for a  classical  particle method and to \cite{Dil99,TK03,TKH09,ZY02} for different developments of the original idea from \cite{GM77}.

In the present work we  extend  a  particle method  developed in \cite{KT15} which is especially adapted to the solution of  
macroscopic mean field equations derived from interacting particle systems.
The method  can be viewed as  a numerical  transition between a DEM-type method for microscopic  interacting particle systems  and a  meshfree particle method for macroscopic equations.
For related approaches we refer to \cite{WTLB12,ZY02}.
In  macroscopic  mean-field models  an interaction term appears, which is  derived from the microscopic interaction term. Usually, this term has the form of a convolution integral. 
The 
particle method  approximates these convolution integrals in an appropriate way not using a microscopically large number of macroscopic particles and thus
allowing the use of 
an underresolved meshfree method.
In the limit for small $R$ a method which is consistent with the asociated limit equation is obtained.
For intermediate values of $R$ one  obtains with this procedure an easily calculated correction term.
In the present work more general particle systems and more general limit equations are considered compared to \cite{KT15} where only porous media type limits have been investigated.
The present investigations include for example terms leading to hyperbolic limit equations and more general applications.

The paper is organized as follows. Section \ref{sec:model} contains a description of the models under investigation. In particular, the macroscopic hydrodynamic  models and their scalar approximations are discussed. Examples range from non-local 1D Lightill Whitham type models to non-local 2D pedestrian flow models.
Section \ref{sec:Approximation} describes the numerical method and in detail the  approximation of the convolution integral. Finally, Section \ref{sec:Numerics} contains the numerical results and a comparison of the methods and their computation times  for the above mentioned examples.

\section{The models} \label{sec:model}

We consider a hierarchy of models ranging from microscopic interacting particle systems and their mean field approximations to hydrodynamic and scalar macroscopic approximations. 

\subsection{Microscopic and mean field equations for interacting particle system}
The starting point for our  model is a general microscopic 
model for $ N$ interacting particles \cite{Jou06}. We define the empirical density as
\begin{align}
\label{empirical}
 \rho^N_R (x) = \frac{1}{N} \sum_{j=1}^N U_R (x - x_j  )
\end{align}
and consider systems of equations of the form
\begin{eqnarray}
\label{micro}
  dx_i &=& v_i dt \\
 dv_i &=& - \left(\beta \left( \rho^N_R (x_i)   \right) 
          +  \alpha\left(  \rho^N_R (x_i)   \right)
          \nabla_x \rho^N_R (x_i)   
             + \gamma v_i \right) dt
            +  A  d W_t \nonumber 
\end{eqnarray}
with $i=1, \ldots, N$,  $(x_i , v_i) \in \R^d \times \R^d$, $d \in \mathbb{N}$. $\alpha: \R^+\rightarrow \R^+,\beta: \R^+\rightarrow \R^d$, $ \gamma, A \ge 0$ and $W_t$ a d-dimensional Brownian motion.
 $U_R$ is a  sufficiently smooth,  repulsive  potential with  support in 
$B_R(0)$ and positive integral which is normalized to $1$, that means $$ \int_{B_R(0)} U_R(x) dx  = 1 . $$
We assume that $U_R$ approximates the delta distribution as $R$ goes to $0$, $U_R (x) = \frac{1}{R^d} U_1(\frac{x}{R})$.  
Possible extensions  are  given in the remarks below.

\begin{remark}
We note that  physically speaking  $U_R$ is used here on the one hand to define the  empirical density,
which is used in the definition of the coefficients $\alpha$ and $ \beta$.
 On the other hand it is used as interaction potential in the term $\nabla_x \rho^N_R$.
 For the following calculations we could as well use different functions $U_R$ at the different places in the equation. 
\end{remark}

 \begin{example}
 \label{quadratic}
As an example for a repulsive interaction potential we choose for $\Vert x \Vert \le R$
\begin{eqnarray}
\label{potsymm}
 U_R =  C_R \left( R -  \Vert x \Vert\right)^2,
 \end{eqnarray}
 where $$C_R = \frac{(d+1)(d+2) }{2 \tau_d R^{d+2}}
  $$ and  $\tau_d$ is the volume of the unit sphere in $\R^d$.
The coefficients are chosen such that 
 $$ \int_{B_R(0)} U_R(x) dx  = 1. $$ 
 We note that non-symmetric potentials could be chosen  as well.
  \end{example}

\begin{remark}
\label{rem1}
In 1D the present formulation includes one sided interaction potentials  $V_R$ approximating the Heaviside function, i.e.
$\partial_x V_R = U_R$. 
An equation of the form 
\begin{eqnarray*}
 dv_i &=& -  \alpha\left(\rho^N_R (x_i ) \right)
          \left( \frac{1}{N}  \sum_{j } \partial_x V_R (x_i - x_j  ) \right) dt 
\end{eqnarray*}
turns into
\begin{eqnarray*}
 dv_i &=& - \beta \left( \rho^N_R (x_i )\right)
           dt
\end{eqnarray*}
with $\beta(\rho) = \rho \alpha(\rho)$.
Compare the systems considered in  \cite{szn,BV05}. 
\end{remark}
 \begin{example}
 \label{unsymm}
  For the potential $V_R$ in the remark above  with $\partial_x V_R =U_R$ one might choose 
   potentials,such that   $U_R$ is  given as in  example 1 or 
    non-anticipating  examples like 
 \begin{eqnarray} 
 \label{potunsymm}
  U_R = 2 C_R \left( R -  \vert x \vert\right)^2,
  \end{eqnarray}
  if $x<0$
  and $0$ otherwise.
   \end{example} 
  
  \begin{remark}
  \label{remmorse}
  For the following considerations, one might also consider more general attractive-repulsive potentials. For example  in 1-D
  the following class of potentials can be considered, compare \cite{LTB09,BV05}:
  \begin{align}
  \label{morse}
  U(x) = G(\vert x \vert) - L \cdot  F G\left(\frac{\vert x \vert}{L}\right),
  \end{align}
  where $G: \R_+ \rightarrow \R_+$ is a monotone  decaying, integrable function, decaying to 0 as $x$ goes to infinity.
  This includes, for example, the Morse potential ($G(x) = exp(-x), x>0$) used among other applications for modelling the swarming behaviour of birds, see \cite{CDP}.
  However, for the following, the potential should lead to a spreading behaviour of the solutions goverened for long times by a diffusive equation,
  in our case a porous media type equation, see \cite{oehl,BV05,LTB09}.
  For the above class of potentials a purely repulsive potential is characterized by $F,L<1$, see \cite{LTB09}.
  However, this is not a necessary condition for the spreading behaviour of the solution. The analysis in \cite{LTB09} shows that one still observes a spreading behaviour, if $F<1$ and $1-F\cdot L^2 >0 $, i.e. also for $L>1$ with $L < 1/\sqrt{F}$.
  We note that in this case 
  $$
  D = \int U(x) dx = 2 \int_0^\infty G(x) dx ( 1-F L^2 ) >0.
  $$
  In $d$ dimensions the corresponding formula for a rotationally symmetric potential $U$ defined by (\ref{morse}) is 
  $$
    D = \int U(x) dx = d \tau_d \int_0^\infty G(r) r^{d-1} dx ( 1-F L^{d+1} ) .
  $$
  This is larger than $0$, if $1-F L^{d+1}>0$.
  \end{remark}

 \begin{remark}
 \label{rem2}
 Velocity dependent interactions $U_R = U_R(x_j-x_i,v_j-v_i)$
 might be included, as well as 
   further explicit dependence of $\alpha$ and $\beta$ on $x_i,v_i$.  
  For example, one could include an exterior potential $V(x)$ via $\beta (x) = \nabla_x V(x)$.
   \end{remark} 
   
  For $N$ going to infinity, one can
  derive in the limit of a large number of particles the associated
  mean field equation \cite{CCR,CDP,BH,neunzert}.
  One obtains  for the distribution function $f=f(x,v,t)$ of the particles
 the mean field equation 
 \begin{eqnarray}
 \label{pde}
   \partial_t f  +  v \cdot \nabla_{x} f =  S f
   + L f 
 \end{eqnarray}
 with force term
 \begin{align*}
 \begin{split}
 S f &= 
     \nabla_v \cdot \left( \beta(  U_R \star \rho ) f\right) 
   +   \nabla_v \cdot \left( \alpha(  U_R \star \rho  )  \nabla_x U_R \star \rho f\right).
     \end{split}
    \end{align*}
 and diffusion term
 \begin{align*}
  \begin{split}
  L f &= \gamma 
     \nabla_{v} \cdot \left(v f +  \frac{A^2}{2 \gamma} \nabla_v f \right).
      \end{split}
     \end{align*}
 Here, the convolution is defined as 
 \begin{align*}
    F_R \star \rho  (x) = 
    \int  F_R(x-y ) \rho(y)   dy   
    \end{align*}
 and the density as 
 \begin{eqnarray*}
 \rho(x,t) := \int  f(x,v,t) d v .
 \end{eqnarray*}
  We  normalize
  \begin{eqnarray*}
  \int \rho(x,t)  d x    =1. 
  \end{eqnarray*}

\subsection{Hydrodynamic and scalar macroscopic models}

Multiplying the mean field equation with $1$ and $v$ and closing the equations by approximating the distribution function with a function having  mean $u$ and variance $ \delta \rho $,   one obtains the continuity and momentum equations
\begin{align}
 \label{cont}
  \partial_t \rho   +   \nabla_{x} \cdot  (\rho u)  &= 0, \\
 \partial_t u + (u \cdot  \nabla_{x}  ) u + \frac{\delta}{\rho} \nabla_x
 \rho &=   - \gamma u  
-  \beta (U_R \star \rho) -   \alpha (U_R \star \rho) \nabla_x U_R \star \rho \nonumber
\end{align}   
with the momentum
 \begin{eqnarray*}
 \rho u(x,t) := \int v  f(x,v,t) d v .
 \end{eqnarray*}
Neglecting time derivatives and inertia terms in the $u$-equation we approximate the velocity $u$ as
$$
u = -\frac{1}{\gamma} \left(  \alpha(U_R \star \rho)\nabla_x U_R \star \rho  + \beta(U_R \star \rho) + \frac{\delta}{\rho} \nabla_x \rho \right)
$$
and obtain  the following scalar equation for the density
 \begin{eqnarray}
   \label{diff1}
   \gamma \partial_t \rho  =   
    \nabla_{x} \cdot \left( \rho \left(\beta(U_R \star \rho)
  +   \alpha(U_R \star \rho) \nabla_x U_R \star \rho  \right) + \delta \nabla_{x}  \rho\right) .
\end{eqnarray}

\subsection{Localized models}

Since  $U_R$ approximates  for small values of $R$  a $\delta$ distribution, 
 one obtains formally
$$
\int \nabla_x U_R(x-y) \rho(y) d y = - \int \nabla_y U_R(x-y) \rho(y) d y = \int U_R(x-y) \nabla_y \rho dy
  \sim   \nabla_x \rho.
   $$ 
   Thus, one obtains  from the hydrodynamic equations the damped isentropic Euler equations with exterior potential
 \begin{align}
 \label{contdiff}
   \partial_t \rho   +   \nabla_{x} \cdot  (\rho u)  &= 0,\\
  \partial_t u + (u \cdot  \nabla_{x}  ) u &=   - \gamma u 
 -     \alpha( \rho) \nabla_x \rho - \beta( \rho) - \frac{\delta}{\rho} \nabla_x \rho. \nonumber
 \end{align}     
 From the   scalar  equations one obtains with the same procedure as before a nonlinear Fokker-Planck equation of the form 
\begin{eqnarray}
  \label{diff2} 
   \gamma \partial_t \rho  =   \nabla_{x} \cdot  (  
      \beta( \rho) \rho )  +    \nabla_{x} \cdot ( (\delta+ \rho \alpha(\rho)) \nabla_x \rho ) . 
\end{eqnarray}
%
%
%
A simple example with a solution which converges to a stationary distribution is given by the following

\begin{example}
	\label{examplebasic}
We consider  for $d=1$ the following coefficients $\alpha(\rho,x) =\epsilon, \beta(\rho,x) =   \partial_x V (x)  + (1-\eps) \rho, \epsilon \in (0, 1]$, $A =  0, \delta =0$.
This leads to the microscopic problem 
\begin{eqnarray*}
  dx_i &=& v_i dt \\
 dv_i &=& - \gamma \left(  \partial_x V (x_i)   +   (1-\eps) \rho^N_R(x_i  )
           +  
          \epsilon  \partial_x   \rho^N_R(x_i  )  
             +  v_i \right) dt 
             \nonumber 
\end{eqnarray*}
and the scalar equation
\begin{eqnarray}
\label{basicnonlocal}
      \partial_t \rho   =   \partial_{x} \cdot ( \epsilon \rho \partial_x U_R \star \rho  +    \rho \partial_x V 
      +  (1-\epsilon) \rho U_R \star \rho ).
 \end{eqnarray}  
The scalar localized macroscopic approximation is 
\begin{eqnarray}
\label{basiclocal}
      \partial_t \rho   =   \partial_{x} \cdot ( \eps \rho \partial_x \rho  +    \rho \partial_x V 
      +  (1-\epsilon) \rho^2 )  .
 \end{eqnarray}  
 with the stationary solution 
 $$
 \rho_{\infty} = \left(C_\epsilon \exp(-\frac{1-\eps}{\eps} x) -    \int^x \frac{1}{\epsilon}\exp(-\frac{1-\eps}{\eps}  (x-x^\prime)) \partial_x V (x^\prime) d x^\prime \right)_+
 $$
 where $C_\epsilon \in \R$ is determined by the normalization $\int \rho_{\infty}(x)  d x =1$.
 $\epsilon =1$ gives the well known solution of the Fokker-Planck porous media equation $\rho_{\infty} = (C- V(x))_+$, see \cite{Aro97}. For $\epsilon $ tending to $0$ one approaches  a  stationary distribution function, which is for a convex potential $V$ given by  $\rho_{\infty} = (-V^\prime (x))_+$ for $ C<x$ and $0$ otherwise.
 $C$ is in both cases  chosen such that the normalization condition is fulfilled.
 \end{example}

We consider further  examples  from traffic and  pedestrian  flow in one- and two space dimensions.

\begin{example}
	\label{exampletraffic}
A simple microscopic traffic model
 is given as follows. 
 We model the acceleration by choosing $\gamma >0$ and $\beta =- \gamma$. Breaking  interactions are  modelled by an interaction potential $V_R$ given by a  smooth version  of the Heaviside function such that $\partial_x V_R   = U_R$, where
 $U_R$ is a smooth version of the   $\delta$-function.
$U_R$ is 
 concentrated in the negative half plane to model the fact that the interaction  is essentially restricted to an interaction with the  predecessors. We call $V_R$  a downwind potential in this case.
 Finally, we choose $\alpha=\gamma$ and  $A  = 0$. 
 We obtain with $x_i,v_i \in \R, i= 1 \cdots N$
\begin{eqnarray*}
  dx_i &=& v_i \\
  dv_i &=& \gamma (1-  v_i  -\rho^N_R (x_i  )) dt \nonumber 
\end{eqnarray*}
The corresponding
mean field hydrodynamic equations are 
\begin{align}
\label{trafficnonlocalhydro}
\partial_t \rho + \partial_x(\rho u )  =0\\
\partial_t u + u \partial_x u + \frac{\delta \gamma}{\rho} \partial_x \rho= \gamma((1- U_R \star \rho) -  u)	\nonumber
\end{align}
and the associated scalar nonlocal viscous Lighthill-Whitham model is
\begin{align}
\label{trafficnonlocal}
\partial_t \rho + \partial_x ((1- U_R \star \rho) \rho) = \delta \partial_{xx} \rho.	
\end{align}
We remark that  nonlocal Lighthill-Whitam model  has been investigated in \cite{GS15,BG15}.
There, for one-sided downwind monotone convolution kernels concentrated in the negative half plane  existence and uniqueness of solutions have been shown as well as 
a maximum principle. In particular, the density is not exceeding the maximal density, which is in this case eqal to $1$ for suitable initial conditions.
The local limits are 
\begin{align}
\label{trafficlocalhydro}
\partial_t \rho + \partial_x(\rho u )  =0\\
\partial_t u + u \partial_x u + \frac{\delta \gamma}{\rho} \partial_x \rho= \gamma(1-  \rho  -  u)
\nonumber
\end{align}
and the Lighthill Whitham equation
\begin{align}
\label{trafficlocal}
\partial_t \rho + \partial_x ((1- \rho) \rho) =\delta \partial_{xx} \rho.	
\end{align}

We note that for the one-sided downwind kernel we obtain an approximation to order $\mathcal{O}(R^2)$ in the following way.
Define
$$
D_R = - \int  x U_R(x) dx = - R \int x U_1 (x) dx >0
$$
using $U_R \star \rho = \rho +  D_R \partial_x  \rho   $ we
obtain  
\begin{align*}
\partial_t \rho + \partial_x ((1-  \rho) \rho) =\delta \partial_{xx} \rho
+  R D_1 \partial_x ( \rho \partial_x  \rho).	
\end{align*}
This is a stable equation for $\delta \ge 0$.
We note that using the potential defined in Example 2 we obtain 
$D_1 =  \frac{1}{4}$.
If instead of using a down-wind interaction potential we consider  symmetric potentials  we obtain $D_1 =0$ and the approximation to order $\mathcal{O}(R^3)$ is
\begin{align}
\label{symmkern}
\partial_t \rho + \partial_x ((1-  \rho) \rho) =\delta \partial_{xx} \rho +
 \frac{R^2}{2} D_2 \partial_x ( \rho \partial_{xx} \rho).	
\end{align}
$$
D_2 =  \int  x^2 U_R(x) dx = R^2 \int x^2 U_1 (x) dx >0.
$$
In this case we need $\delta >0$ in order to obtain a stable convergence of the solutions of equation (\ref{symmkern}) as $R$ goes to $0$, compare \cite{scho}. 
A numerical investigation of non-local Lighthill-Whitham type equations with symmetric or upwind potentials
can be found  in \cite{BG15}. In this case there is no maximum principle and the  density might exceed the 'maximal density'.
\end{example}
%
%
%
\begin{example}
For the 2D case we consider a model for pedestrian flow.
For $x_i,v_i \in \R^2, i=1, \ldots, N$ and $W(\rho) = u_{max}(1-\frac{\rho}{\rho_{max}})$  we
 consider the microscopic model
\begin{eqnarray*}
  \frac{dx_i}{dt} &=& v_i dt\\
 dv_i &=& - \gamma W \left( \rho_R^N (x_i  ) \right) \hat e (x_i) dt -  \gamma v_i dt-  
       \alpha  \nabla_x \rho^N_{R} (x_i  )  dt + A dW_t
            \nonumber 
\end{eqnarray*}
with symmetric interaction kernel $U_R$ and
with
$$
\hat e (x) = \frac{\nabla \phi (x)}{\vert \nabla \phi (x) \vert}
$$
where
$$
\vert \nabla \phi (x_i) \vert = \frac{1}{W \left(\rho^N_R (x_i  ) \right)}.
$$
The hydrodynamic mean field limit is
\begin{align}
\label{pednonlocalhydro}
\partial_t \rho + \nabla_x(\rho u )  =0\\
\partial_t u + u \nabla_x u + \frac{\delta}{\rho} \nabla_x \rho= - \gamma \left(W ( U_R \star \rho) \hat e(x) -  u	\right)
- \alpha \nabla_x U_{{R}} \star \rho  \nonumber \end{align}
and the scalar limit is
\begin{align}
\label{pednonlocalscalar}
\partial_t \rho -  \nabla_x \cdot (W( U_R \star \rho) \hat e (x) \rho) = \frac{\delta }{\gamma}
\Delta \rho  + 
 \nabla \cdot \left( \frac{\alpha}{\gamma} \rho \nabla_x U_{{R}} \star \rho  \right) 	
\end{align} 
where
$$
\vert \nabla \phi  \vert = \frac{1}{W(U_R \star \rho)}
$$
The local limits are
\begin{align}
\label{pedlocalhydro}
\partial_t \rho + \nabla_x(\rho u )  =0\\
\partial_t u + u \nabla_x u +\frac{\delta}{\rho} \nabla_x \rho = - \gamma \left(W(  \rho)\hat e  -  u\right) - \alpha \nabla_x \rho \nonumber
\end{align}
and
\begin{align}
\label{pedlocalscalar}
\partial_t \rho - \nabla_x \cdot (W( \rho) \hat e \rho) =
 \nabla \cdot \left( \frac{\alpha \rho +\delta }{\gamma} \nabla_x  \rho  \right)
\end{align}
together with
$$
\vert \nabla \phi  \vert = \frac{1}{W(\rho)}.
$$
which is a viscous form of the Hughes model.
For a detailed  modelling of the interactions between pedestrians we refer for example to \cite{Deg13}.
For rigorous results on the Hughes model and the approximation of the model via particle systems we refer to \cite{di,hughes1,di2}.
\end{example}

In the following, our goal will be to develop a meshfree particle method
for equations (\ref{cont}), (\ref{diff1}) for different ranges of parameters
and, in particular, for the limit equations
(\ref{contdiff}), (\ref{diff2}).

\section{Numerical method}

For  a general description of the particle method used here and further references on the subject, we refer to  \cite{KT15}. Here, we concentrate on the approximation of the interaction term.
We follow the approach in \cite{KT15}. However, additionally to \cite{KT15} we have to treat here the terms
leading to hyperbolic limit equations in a suitable way.

\subsection{Approximation of the interaction term using particle methods}

\label{sec:Approximation}
The key point of our method is to approximate the integrals
\begin{align}
\label{int1}
 F_R \star \rho (x) = \int  F_R (x-y) \rho (y) dy
\end{align}
with $F_R= U_R$ and $F_R = \nabla_x U_R$.
We denote the Voronoi cell around particle $i$ given by the particle locations 
of the other particles by $V_i$. $\vert V_i \vert$ denotes the volume of this cell.
Then, a  naive or microscopic approximation of the integral terms  would be to use
\begin{align}
\label{simple}
F_R \star \rho (x_i) \sim  \sum_{j=1, j \neq i}^N \rho_j \vert V_j \vert F_R(x_i-x_j)  , 
\end{align}
where $N$ is now the number of macroscopic particles used in the particle method.
For small values of $R$ this results in the following  problem. We consider a situation where  the method is underresolved, that means where  one uses a small number of macroscopic particles
in contrast to the large number of microscopic particles described by the macroscopic equations.
Then, the value evaluated from (\ref{simple}) will be zero due to the large distances between the macroscopic particles and the relatively small value of $R$.
However,  the actual 
value of the integral will not be zero even for very small $R$
due to the corresponding 'infinite' number of microscopic particles described by the macroscopic equations. 

We resolve this problem using a higher order approximation of the integral.
This yields in the limit for  small $R$ a method for the limiting nonlinear Fokker-Planck equations,  even if the number of macroscopic particles is still small.
We use an approximation of the density given by 
\begin{align}
\rho(y) = \sum_{j=1}^N \left[\rho_j + \sigma_j \cdot (y-x_j)\right] \chi_{V_j} (y),
\end{align}
where $V_j$ denotes the Voronoi cell associated to particle/mesh point $x_j$
and $\chi$ denotes the characteristic function.
The approximation $\sigma_j$  of the first derivative 
$\nabla_x \rho (x_j)$  is determined via a  least squares approximation  using the neighbouring points.
Then, we obtain for the integral (\ref{int1})
\begin{align}
\label{int}
& 
F_R \star \rho (x_i) \\
&\sim   \sum_{j=1}^N \left[\rho_j \int_{V_j \cap B_R(x_i)} F_R(x_i-y) dy   + \sigma_j \cdot \int_{V_j \cap B_R(x_i)} (y-x_j) F_R(x_i-y)  dy \right] .\nonumber
\end{align}
We approximate first the integral around the center point $x_i$.
 To approximate the integral 
 we distinguish between the cases $\vert V_i \vert > \vert B_R \vert$
 and $\vert V_i \vert < \vert B_R \vert$.
 For $\vert V_i \vert > \vert B_R \vert$, we proceed as follows.
  \begin{align*}
 \int_{V_i \cap B_R(x_i)} F_R(x_i-y) dy  \sim
 \int_{ B_R(x_i)} F_R(x_i-y) dy =1
 \end{align*}
 This is equal to 1 for $F_R = U_R $ due to the normalization of the potential.
 For  $F_R = \nabla_x U_R $ the expression is equal to 0, since
  \begin{align*}
  \int_{ B_R(x_i)} \nabla_x U_R (x_i-y) dy  =
  \int_{B_{R}(0)} \nabla_x U_R (y) dy 
  = \int_{ \partial B_{R}(0)} \frac{y}{\Vert y \Vert}  U_R(y) d S(y) =0.
  \end{align*} 
 Moreover,
\begin{align*}
 \int_{V_i \cap B_R(x_i)} \sigma_i \cdot (y-x_i)  F_R(x_i-y) dy \sim 
\int_{ B_R(x_i)} \sigma_i \cdot (y-x_i)  F_R(x_i-y) dy  \\
= -  \int_{ B_R(0)} \sigma_i \cdot y F_R(y) dy  .
\end{align*} 
This expression is for $F_R = U_R$ equal to $-  \sigma_i \cdot \mu$, where $\mu = \int_{ B_R(0)}  y U_R(y) dy $ is the mean value of $U_R$. For $F_R = \nabla_x U_R$ we obtain
 \begin{align*}
 -  \int_{ B_R(0)} \sigma_i \cdot y  \nabla_x U_R (y) dy 
 =     \int_{ B_R(0)} U_R(y)  dy \sigma_i
 =     \sigma_i. 
\end{align*}
If $\vert V_i \vert < \vert B_R \vert$ then we first compute
$R_{V_i}$ such that $\vert B_{R_{V_i}} \vert  = \vert V_i \vert$. Then,
 \begin{align*}
 \int_{V_i \cap B_R(x_i)} F_R(x_i-y) dy  \sim
 \int_{V_i } F_R(x_i-y) dy  \\\sim
 \int_{B_{R_{V_i}}(x_i)} F_R(x_i-y) dy = \int_{B_{R_{V_i}}(0)} F_R(y) dy \nonumber
 \end{align*}
Thus we obtain $\int_{B_{R_{V_i}}(0)} U_R(y) dy$ for $F_R = U_R$. 
For $F_R = \nabla_x U_R$ we get 
\begin{align*}
 \int_{ \partial B_{R_{V_i}}(0)} \frac{y}{\Vert y \Vert}  U_R(y) d S(y)
\end{align*} 
Finally, we compute 
\begin{align*}
 \int_{V_i \cap B_R(x_i)} \sigma_i \cdot (y-x_i)  F_R(x_i-y) dy \sim 
 \int_{V_i } \sigma_i \cdot (y-x_i)  F_R(x_i-y) dy \\
 \sim 
- \sigma_i \cdot \int_{B_{R_{V_i}}(0)}  y  F_R(y) dy.\nonumber
\end{align*} 
This is equal to $- \sigma_i \cdot \int_{B_{R_{V_i}}(0)}  y  F_R(y) dy$
for $F_R = U_R$. 
For $F_R = \nabla_x U_R$ we get 
 \begin{align*}
 &\sigma_i \int_{ B_{R_{V_i}}(0)}  U_R(y) dy - \int_{ \partial B_{R_{V_i}}(0)} \sigma_i \cdot y \frac{y}{\Vert y \Vert}  U_R(y) d S(y) . 
 \end{align*}
The integrals over the Voronoi cells  with the points $x_j, j \neq i$ as centerpoints are approximated in the following way.
A simple second order approximation is given by 
the midpoint rule
 \begin{align*}
\int_{V_j \cap B_R(x_i)} F_R(x_i-y) dy \sim
\vert V_j \vert  F_R(x_i-x_j)
\end{align*}
and
\begin{align*}
 \int_{V_j \cap B_R(x_i)} (y-x_j) \otimes F_R(x_i-y) dy \sim 0.
 \end{align*}
Altogether, one obtains 
\begin{align}
\label{approxu}
U_R \star \rho (x_i) \sim  \sum_{ j \neq i} \Big(\rho_j \vert V_j \vert U_R(x_i-x_j) \Big) +  \alpha_i \rho_i + \alpha^\prime_i \cdot \sigma_i  
\end{align}
and
\begin{align}
\label{approxf}
\nabla_x U_R \star \rho (x_i) \sim  \sum_{ j \neq i} \Big(\rho_j \vert V_j \vert \nabla_x U_R(x_i-x_j) \Big) + \beta_i \rho_i +  \beta^\prime_i  \sigma_i
\end{align}
with the correction factors
\begin{align}
\label{alphaformula}
\alpha_i  =   \int_{B_{R_{V_i}}(0)} U_R(y) dy \\
\alpha^\prime_i  =   -  \int_{B_{R_{V_i}}(0)}  y  U_R(y) dy
\end{align}
\begin{align}
\label{betaformula}
\beta_i  =  \int_{ \partial B_{R_{V_i}}(0)} \frac{y}{\Vert y \Vert}  U_R(y) d S(y) \\
\beta^\prime_i  = \int_{ B_{R_{V_i}}(0)}  U_R(y) dy I
- \int_{ \partial B_{R_{V_i}}(0)}  y \otimes y \frac{1}{\Vert y \Vert}  U_R(y) d S(y) 
\end{align}
for $\vert V_i \vert < \vert B_R \vert$ and 
\begin{align}
\label{alphaformula2}
\alpha_i  =   1 \\
\alpha^\prime_i  =  -  \int_{B_{R}(0)}  y  U_R(y) dy
\end{align}
\begin{align}
\label{betaformula2}
\beta_i  =  0 \\
\beta^\prime_i  =  I 
\end{align}
for $\vert V_i \vert < \vert B_R \vert$.
We note that in order to obtain a stable approximation we have to guarantee that 
$\alpha_i^\prime$ and $\beta_i\prime$ are non-negative.
Alltogether, we have the  following
algorithm
\begin{algorithm}
\begin{enumerate}
\item
For each $x_i$ check for particles inside the interaction-region  $B_R(x_i)$.
\item
Compute for all particles inside this region the size of the Voronoi cells
 $
\vert V_j \vert , j = 1, \ldots, N.
$
\item
Use the  neighbouring particles inside the neighbourhood region of radius $h$
(not only the interaction region) to compute an approximation 
$\sigma_i$ of the first derivative $\nabla_x \rho(x_i)$.
\item
Compute $
F_R \star \rho (x_i) $ according to the above formulas.
\end{enumerate}
\end{algorithm}

\begin{remark}
For $\vert V_j \vert << \vert B_R \vert$ the approximation
behaves like the microscopic interaction approximation, since $\alpha_i, \beta_i$ go to $0$.
For $\vert V_j \vert >> \vert B_R \vert$, i.e. the underresolved situation, we have
an approximation which behaves like a solution method for the macroscopic equations
since $\alpha_i , \beta_i$ go to $1$ and the other terms vanish.
\end{remark}

\begin{remark}
More accurate approximations of the integrals are possible at the expense 
of a more complicated approximation, compare \cite{KT15}.
\end{remark}

\begin{remark}
 	In a situation as in the above remark  with $\beta \ne 0$  we need an upwind procedure to stabilize the numerical approximation. A  first order upwind procedure  amounts  to adding  numerical diffusion proportional to 
 	$
 	\vert V_i \vert     .
 	$
 	In comparison, the above factor $\alpha_i^\prime$, describing the physical diffusion in case of an unsymmetric potential, is  proportional to  $R$. In  case $\vert V_i \vert \sim R$, the two diffusion coefficients 
 	are of the same order.
 	Thus, the physical  diffusion  is of the same order as the numerical diffusion. In order to capture the effects of the physical diffusion a higher order upwinding procedure would be necessary to reduce the numerical diffusion.
 \end{remark}

\subsection{Special cases}

 The above formulas  give for radially symmetric potentials 
 \begin{align*} 
   \alpha_i =    d  \tau_d \int_0^{R_{V_i}} U_R(r) r^{d-1} dr  \rho_i 
 \end{align*}
 and 
\begin{align*} 
 \beta_i ^\prime = d  \tau_d \int_0^{R_{V_i}} U_R(r) r^{d-1} dr  \sigma_i -
    \tau_d   R_{V_i}^d  U_R (R_{V_i}) \sigma_i  ,
\end{align*}
where $\tau_d$ is the volume of the unit sphere in $\R^d$.
$\alpha_i^\prime$ and $\beta_i$  are $0$ in the radially symmetric  case.
Rewriting gives
 \begin{align}
      \alpha_i =   A \left(\frac{R_{V_i}}{R}\right)
      \end{align}
      with
      \begin{align}
                  A(z) =    d \tau_d \int_0^z U_1(y) y^{d-1} d y 
                  \end{align}
       for $z \in [0,1]$ and   $A (z) =1, z \ge 1$.
       Moreover, 
       \begin{align}
           \beta_i^\prime =   B^\prime \left(\frac{R_{V_i}}{R}\right)
           \end{align}
           with
           \begin{align}
                             B^\prime(z) =    d \tau_d \int_0^z U_1(y) y^{d-1} d y - \tau_d z^d U_1(z)
                             \end{align}
        for $z \in [0,1]$ and   $B^\prime (z) =1, z \ge 1$.

\begin{remark}
We note that $\beta_i^\prime$  is positive, if
the potential fulfills
$$
d \int_0^z U_1 (y) y^{d-1} d y - z^d U_1(z) \ge 0, z \in [0,1]
$$
or
$$
\int_0^z y^d U_1^\prime (y) dy  \le 0, z \in [0,1].
$$
This is the case for any repulsive potential. Also potentials like the Morse potential fulfilling the condition in Remark \ref{remmorse} fulfill the condition, since
\begin{align*}
\int_0^z y^d U^\prime (y) dy & = 
\int_0^z y^d \left(G (y) -F LG (y/L)\right)^\prime dy\\
= \int_0^z y^d \left(G^\prime (y) -F G^\prime (y/L)\right) dy
&= 
\int_0^z y^d G^\prime (y) dy
- F L^{d+1} \int_0^{z/L} y^d G^\prime (y) dy .
\end{align*}
This expression is negative for all $z \in [0,\infty]$, if  $F<1$ and 
$1-F L^{d+1} >0 .$
\end{remark}

 	For  the quadratic potential  $U_R$ described in the first section in Remark \ref{quadratic} the above formulas
 	  give
 	the following expressions for $\alpha_i,\beta_i^\prime$:
 	\begin{align}
 	\alpha_i =   A \left(\frac{R_{V_i}}{R}\right)
 	\end{align}
 	with
 	\begin{align}
 	A(z) =  \frac{1}{2}(d+1)(d+2)z^{d} - d (d+2)z^{d+1}+  \frac{1}{2} d (d+1) z^{d+2}  
 	\end{align}
 	for $z \in [0,1]$ and   $A (z) =1, z \ge 1$.
 	Moreover, 
 	\begin{align}
 	\beta_i^\prime =   B^\prime \left(\frac{R_{V_i}}{R}\right)
 	\end{align}
 	with
 	\begin{align}
 	B^\prime(z) =     (d+2) z^{d+1}-   (d+1) z^{d+2}
 	\end{align}
 	for $z \in [0,1]$ and   $B^\prime (z) =1, z \ge 1$.     
 	We note that $B^\prime \ge 0$.         
 	Moreover, $\alpha_i^\prime= \beta_i=0$.
 	
 	In 1D this gives 
 	$$U_R = \frac{3 }{2 R^3} \left(  R - \vert x  \vert 
 	\right)^2, \vert x \vert  \le R$$
 	and the following expressions
 	\begin{align} 
 	\alpha_i =   A  \left(\frac{\vert V_i \vert}{ 2 R}\right), \;
 	A(z) = 3 z -  3 z^2+ z^3 
 	\end{align}
 	and 
 	\begin{align} 
 	\beta_i^\prime =  B^\prime \left(\frac{\vert V_i \vert}{ 2 R}\right),\;
 	B^\prime(z)  =  3  z^2 - 2  z^3   
 	\end{align}
 	for $z \in [0,1]$.
 	In 2-D we have  
 	$$
 	U_R = \frac{6 }{ \pi R^4} \left( R  -  \Vert x \Vert \right)^2,
 	\Vert x \Vert<  R
 	$$
 	and the following expressions
 	\begin{align} 
 	\alpha_i =   A  \left(\sqrt{\frac{\vert V_i \vert}{ \pi R^2}}\right), \;
 	A(z)  =   6z^{2} - 8 z^{3}+ 3 z^{4}
 	\end{align}
 	and 
 	\begin{align} 
 	\beta_i^\prime =  B^\prime \left(\sqrt{\frac{\vert V_i \vert}{ \pi R^2}}\right), \;
 	B^\prime(z) =   4  z^3- 3  z^4.
 	\end{align}
 	for $z \in [0,1]$.

 \begin{remark}	
 Considering the  1-D case,
 	the potential  in Example \ref{unsymm}, i.e.
 	$$U_R = \frac{3 }{ R^3} \left(  R - \vert x  \vert 
 	\right)^2$$ for $ -R <x< 0$ and $0$  otherwise,
 	leads to the expressions
 	\begin{align} 
 	\alpha_i =   A(\frac{\vert V_i \vert}{2R})
 	\end{align}
 	with
 	$$
 	A(z) =
 	3  z - 3 z^2 + z^3, z \in [0,1] $$ and $A(z)=1$ for $ z \ge 1$   
 	and
 	\begin{align} 
 	\alpha_i^\prime  =   \frac{R}{4}  A^\prime(\frac{\vert V_i \vert}{2R})
 	\end{align}
 	with
 	$$
 	A^\prime(z) =
 	6  z^2 -  8 z^3 + 3 z^4 , z \in [0,1] $$ and $A^\prime (z)=1$ for $z \ge 1 $. 
 	We note that  $A^\prime \ge 0$. 	
 \end{remark}

\section{Numerical results}

\label{sec:Numerics}

 In this section we present a series of numerical experiments for 
 the hydrodynamic (\ref{cont}) and scalar  equations  (\ref{diff1}) and  their localized  approximations (\ref{contdiff})   and (\ref{diff2}).  
 In particular, we will investigate and compare the different  schemes for the non-local equations for situations near the local limit.
One observes the following: if a well resolved situation with a large number of particles compared to the interaction radius is considered,
then the results of the  multi-scale and the naive or microscopic
approximation  coincide.
In this case the multiscale method behaves similiar to  a 
microscopic (DEM-type) simulation.  If, however, 
a strongly underresolved situation is considered, i.e. the  number of grid particles  is small compared to the interaction radius, then the 
 numerical  solution using the microscopic approximation of the integral  deviates strongly from the correct solution. On the contrary, using the
 multi-scale method with the correction factors introduced above, a 
 good approximation of the correct  solution is obtained.
 In this case the multi-scale method is essentially a meshfree
 numerical method for macroscopic equations.
Since all methods require approximately the same amount of computation time per grid-particle,
the above observations can be rephrased as follows:
 for situations with relatively small interaction radius we obtain, using the multiscale method, a  reduction in computation time by several orders of magnitude. On the other hand,
for relatively large interaction radius, the computation times of the  microscopic approximation
and of  the multiscale method are similar.

 \subsection{Numerical results in 1-D}

 \subsubsection{Test case (Example 3)}
 We study 1-D movement of the particles under the influence of a  confining potential. 
 The equations are described in Example \ref{examplebasic}. More precisely, we look at equations (\ref{basicnonlocal}) and its limit equation  (\ref{basiclocal}). The interaction potential is given by (\ref{potsymm}). 
 We consider the case of small $\epsilon << 1$. The porous media case $\epsilon =1$ has been treated in Ref. \cite{KT15}.
 We choose $V(x) = \vert x \vert^2/2$ and  $\gamma = 1$.  We choose  different values for the interaction radius  $R$.  
     For this example the particles are assumed to be distributed on
     a fixed equidistant grid, i.e. we use here an Eulerian approach. Situations with arbitrary particle locations and a Lagrangian approach are considered in the following subsections.        
 The initial density is given by  $\rho(x,0) = 1$ for $x \in [-0.5,0.5]$ and $0$ elsewhere.  
  As time proceeds the density   converges 
  to a  stationary solution depending on the value of $R$. For the localized limit equation this stationary solution is for small values of $\epsilon$ approximated by  $\rho_{\infty} = -x $ for $ - \sqrt{2}<x<0$ and $\rho_{\infty} =0$ otherwise, such that $\int \rho_{\infty} (x) dx =1$. In order to obtain the steady state solution we have solved the equation until $t=8$.

  \begin{figure}
     \captionsetup[subfigure]{margin=5pt} 
      \subfloat[$R= 0.002 $]{
        \includegraphics[keepaspectratio=true, width=.5\textwidth]{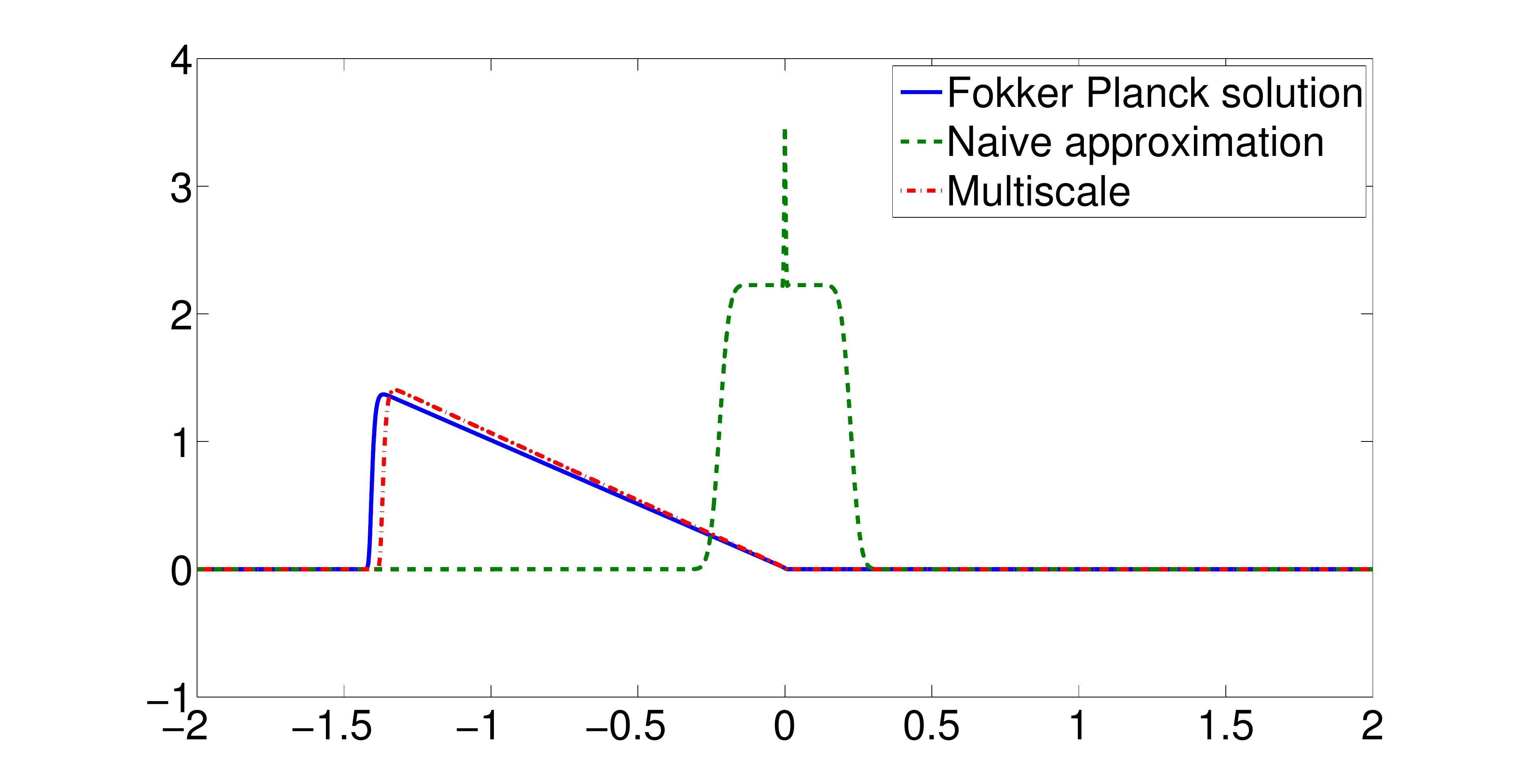}
     } \subfloat[ $R = 0.01$]{
     \includegraphics[keepaspectratio=true, width=.5\textwidth]{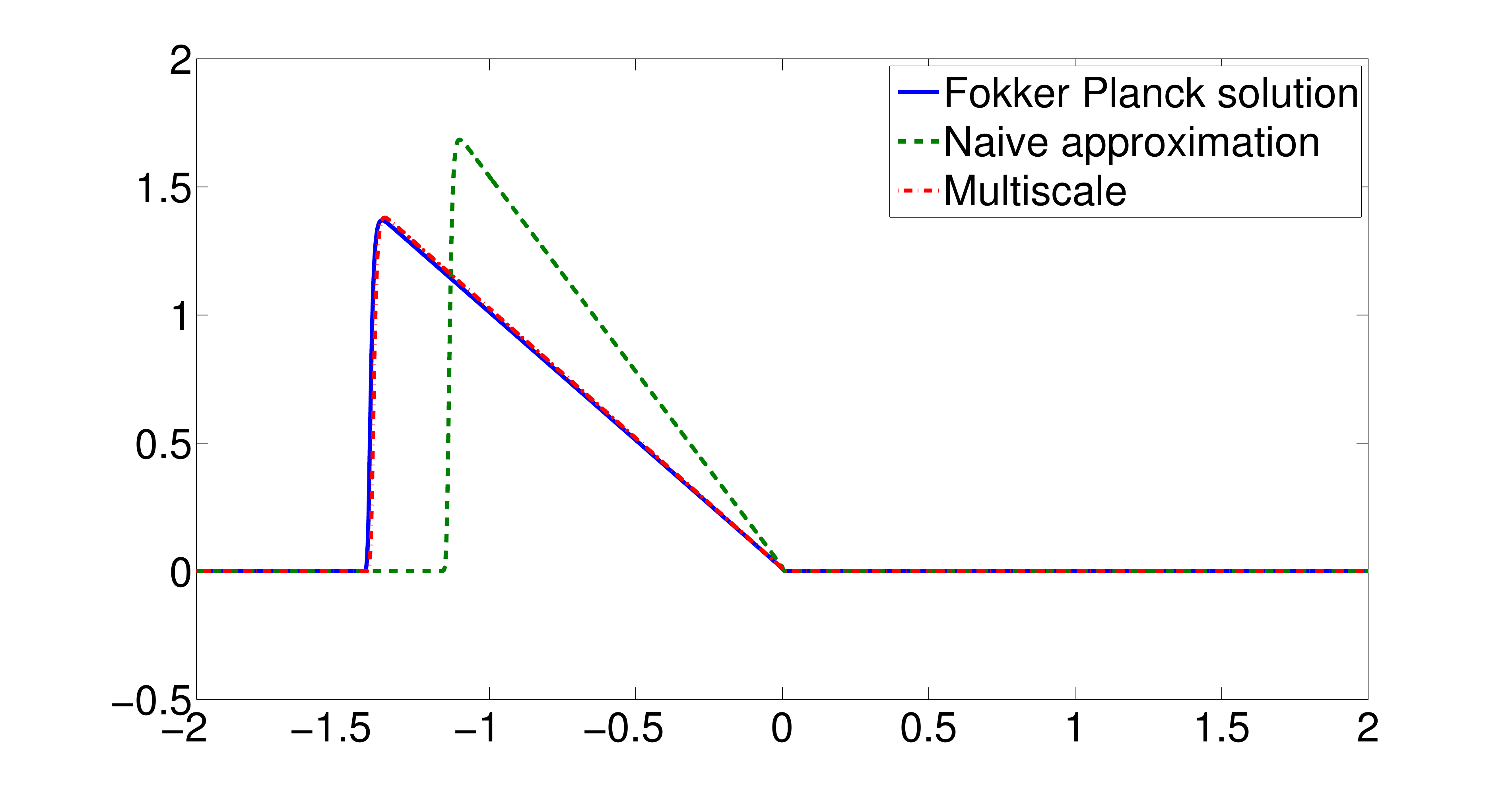}
     } 
      \captionsetup[subfigure]{margin=0pt} 
      \subfloat[ $R =   0.02$]{
      	\includegraphics[keepaspectratio=true, width=.5\textwidth]{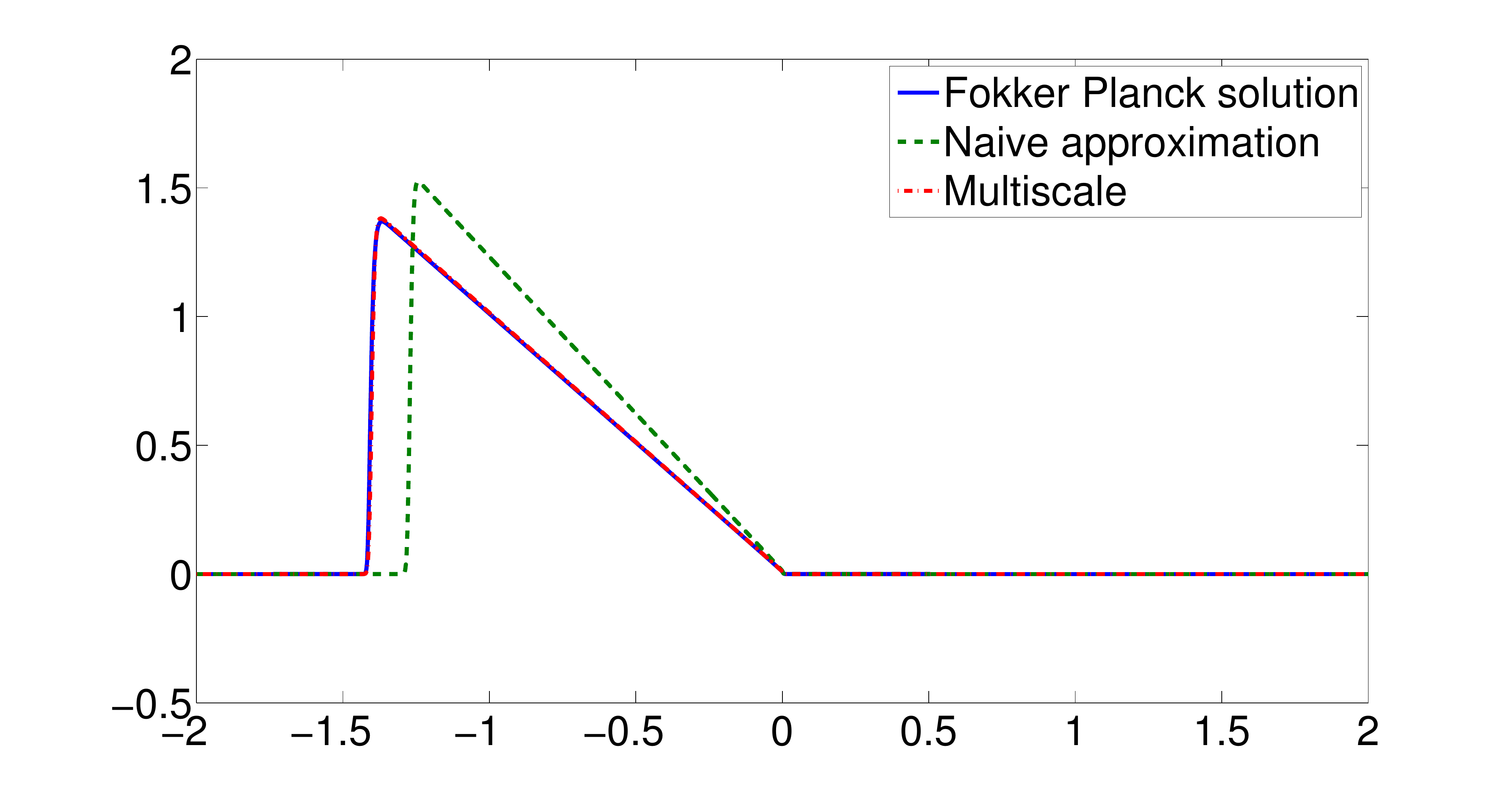}
      }
      \subfloat[$R = 0.04 $]{
      	\includegraphics[keepaspectratio=true, width=.5\textwidth]{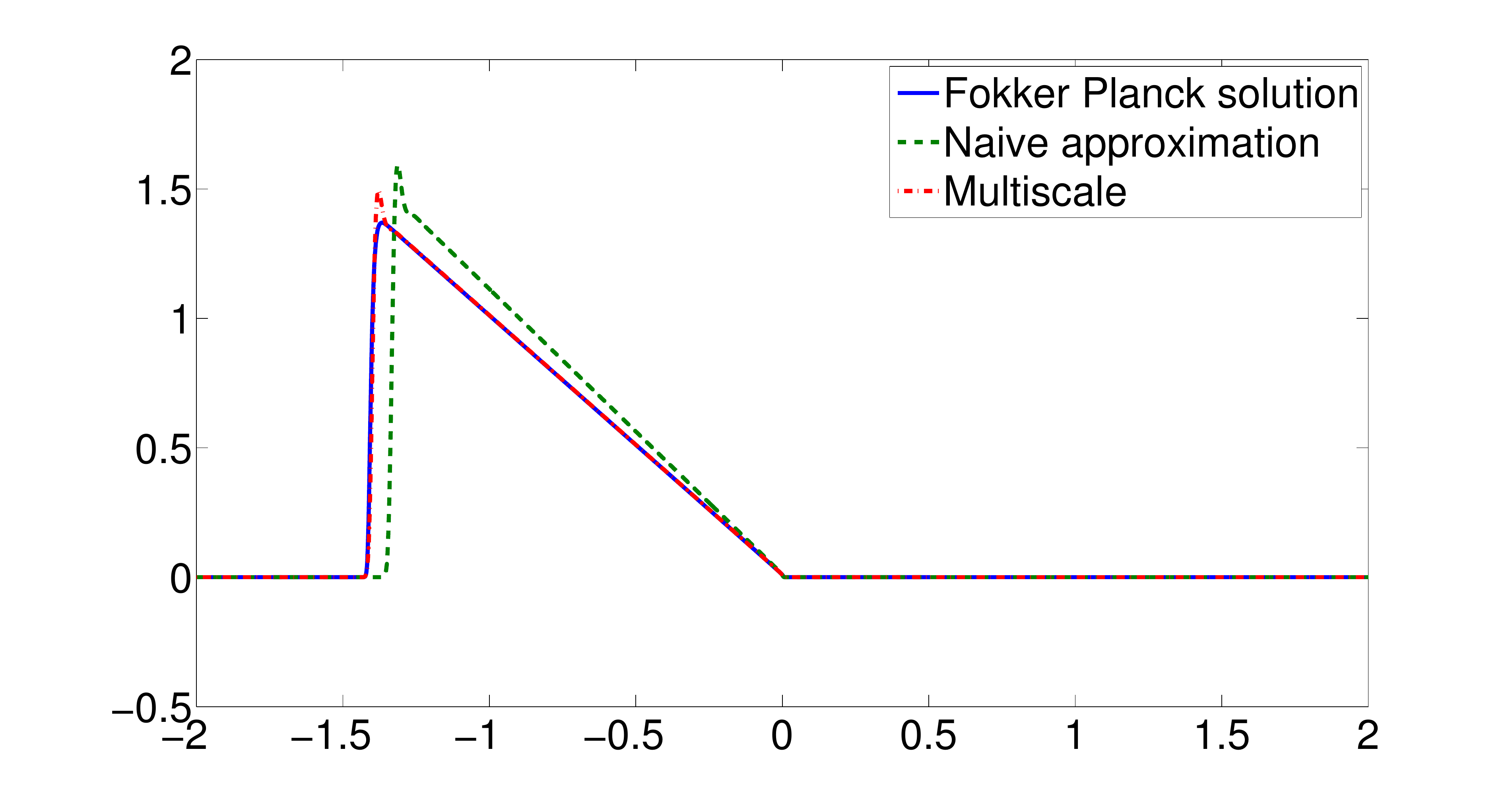}
      }   \\
      \subfloat[ $R = 0.2$]{
      	\includegraphics[keepaspectratio=true, width=.5\textwidth]{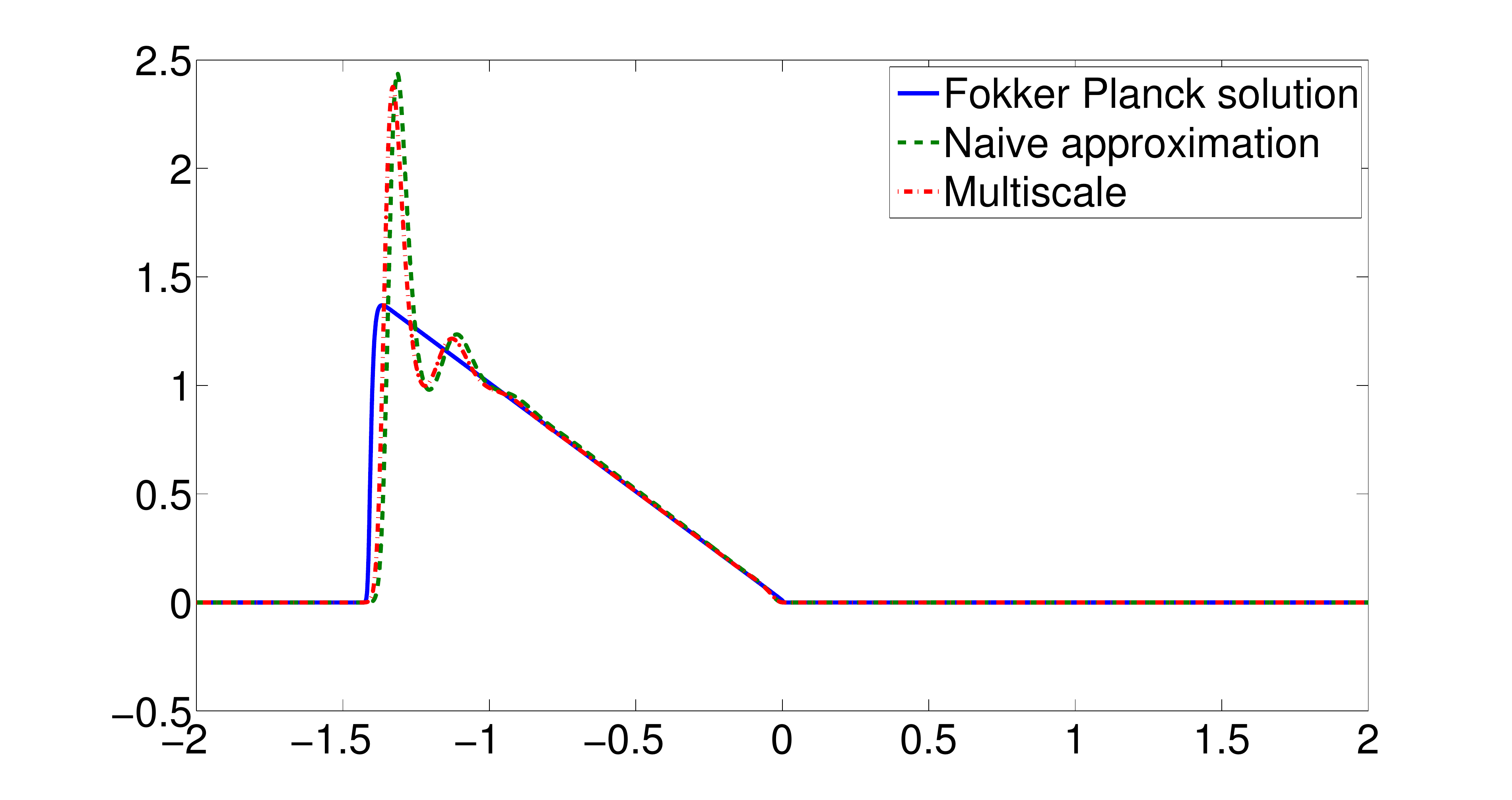}}   
      \subfloat[ $R = 1.0$]{
      	\includegraphics[keepaspectratio=true, width=.5\textwidth]{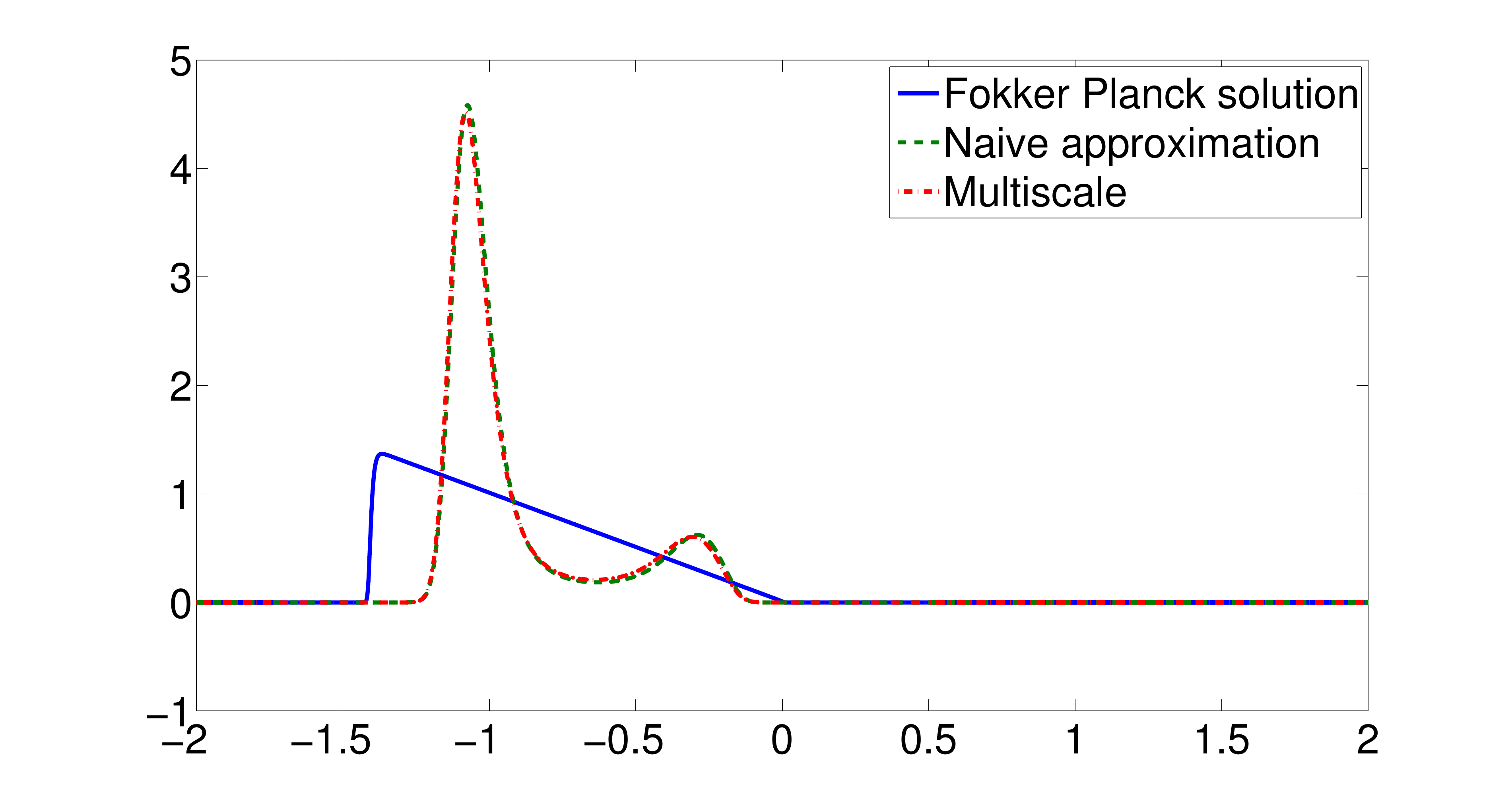}}  
     \captionsetup{margin=20pt}
     \caption{Stationary solution for the  local Fokker-Planck model (\ref{basiclocal}) and the non-local model (\ref{basicnonlocal}) with microscopic and multi-scale approximations.  The interaction radius $R$ is ranging from  $R=0.002$ to  $R=1.0$.    The number of particles is $N=1600$.}
     \label{fig1}
     \end{figure}

    \begin{table}
 \begin{center}
 \begin{tabular}{|r|r|r|r|r|}
  \hline
 $\# $ Particles  & Fokker-Plank &naive  &    multi-scale    & CPU time \\ 
 $ $  & error & error  &    error  &seconds  \\ \hline 
   $200$ & 0.29  & $5.99$    & $0.51$ &$ 26$   \\
 $400$ & 0.21   & $1.14$   & $0.33$  &$ 51$    \\ 
 $800$  & 0.15  & $0.77$  & $0.21$  &$ 104$  \\
 $1600$ & 0.07 & $0.55$   & $0.10$  &$ 216$    \\   \hline
 \end{tabular}
 \caption{$\mathcal{L}^2$-errors and CPU time for Example 3 with $R=0.02$.}
 \label{table1}
 \end{center}
 \end{table}

      \begin{figure}
          \centering
             \includegraphics[keepaspectratio=true, width=.5\textwidth]{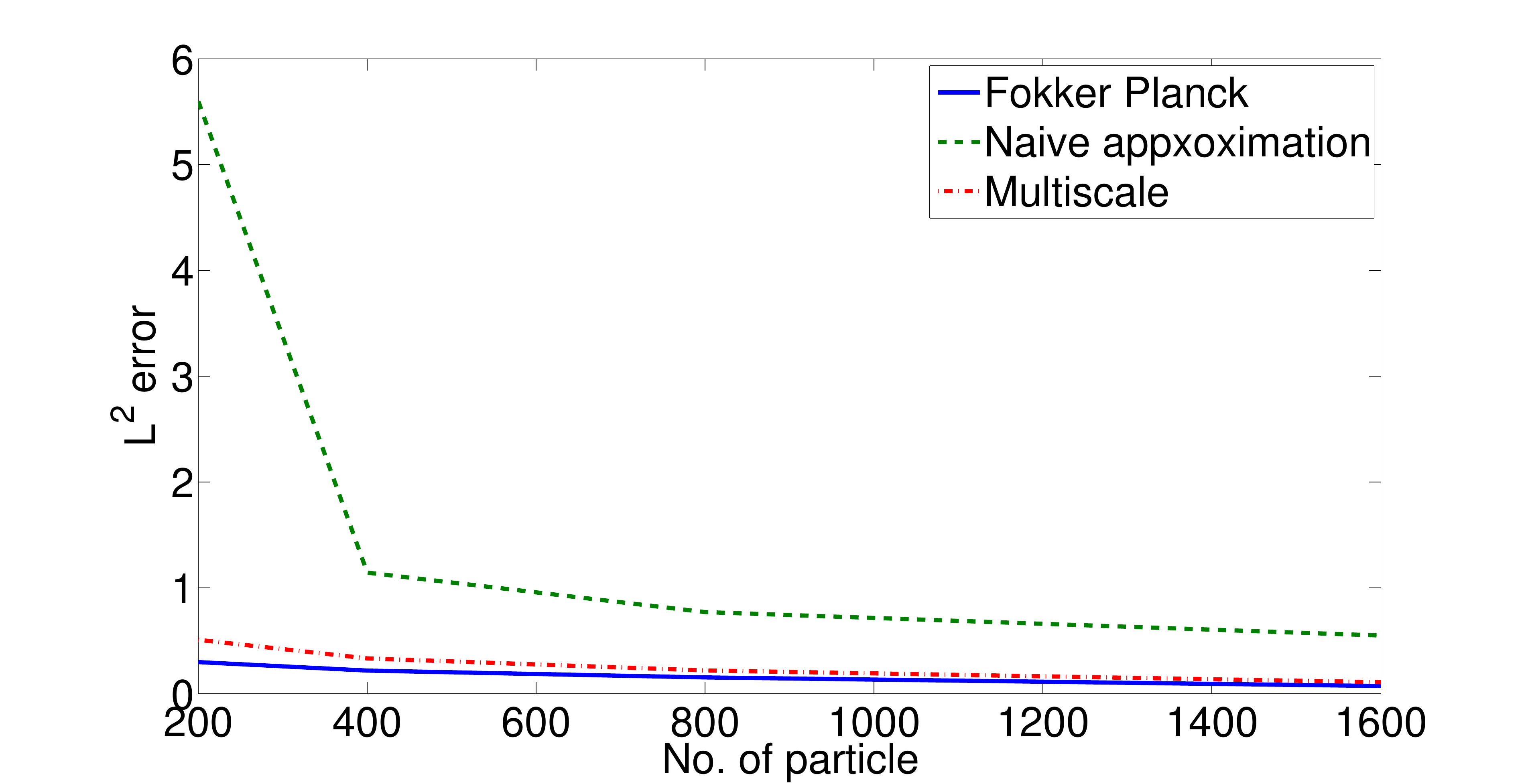}
             \caption{$\mathcal{L}^2$-error plot versus number of particles  for the solutions at steady state  computed from   (\ref{basiclocal}) and  (\ref{basicnonlocal}) 
                                   with $R= 0.02$ using the microscopic and the multi-scale  approximation.}
             \label{error1}
             \end{figure}


 In Figure \ref{fig1} we plot the solutions 
 for a fixed number of particles $N=1600$  and  different values of  $R$,
 considering in this way 
  the  different algorithms   in the well-resolved, in intermediate cases  and in the underresolved case. 
We plot  the limiting Fokker-Planck solution and the microscopic  and multi-scale solutions.
 In the well resolved case, i.e. here for larger values of $R$, the microscopic  approximation of the integral and the multi-scale solution 
 give similiar results deviating from the solution of the limiting Fokker-Planck equation. For the under-resolved case, for  smaller values of $R$, the microscopic approximation  does not give the correct results as discussed above, in contrast to the  the multi-scale method.
 In particular, one observes that using the microscopic discretization for an underresolved situation, the
    numerical  evolution terminates before the stationary state is reached, if  the compact supports of the interaction kernels of the finite number of particles do not overlap any more.

  Finally, we investigate the convergence rate of the method numerically. 
  We consider  a  situation near the localized limit with $R=0.02$, where the limit equation gives the correct solution up to order $10^{-2}$ and is used as the reference solution.
  We compute the error of the particle methods for   the  scalar model (\ref{basicnonlocal}) comparing the numerical results to the numerically determined stationary solution. The scalar problem is solved with microscopic and multi-scale approximation of the interaction term. In Figure \ref{error1}  and Table \ref{table1} the $\mathcal{L}^2$-errors  versus the number of particles are plotted.       
  One observes the deterioration of the   method using the microscopic approximation of the interaction term for smaller numbers of particles.
  The multi-scale method is able to treat all ranges    with a similiar  accuracy.      
  We remark that  the computation time for the different methods is approximately the same for the same number of particles. Looking at Table \ref{table1} we observe that  the multiscale method with $N=200$ particles
yields for the present example the same error as the microscopic method with $N=1600$ particles.
Thus,  the computation time for the 
naive, microscopic method is approximately an order of magnitude larger than the time for the multi-scale method.

%
%
%

 \subsubsection{Non-local traffic flow (Example 4)}

\begin{figure}
    \centering
    \captionsetup[subfigure]{margin=0pt} 
     \subfloat[ $R  = 0.002$]{
        \includegraphics[keepaspectratio=true, width=.5\textwidth]{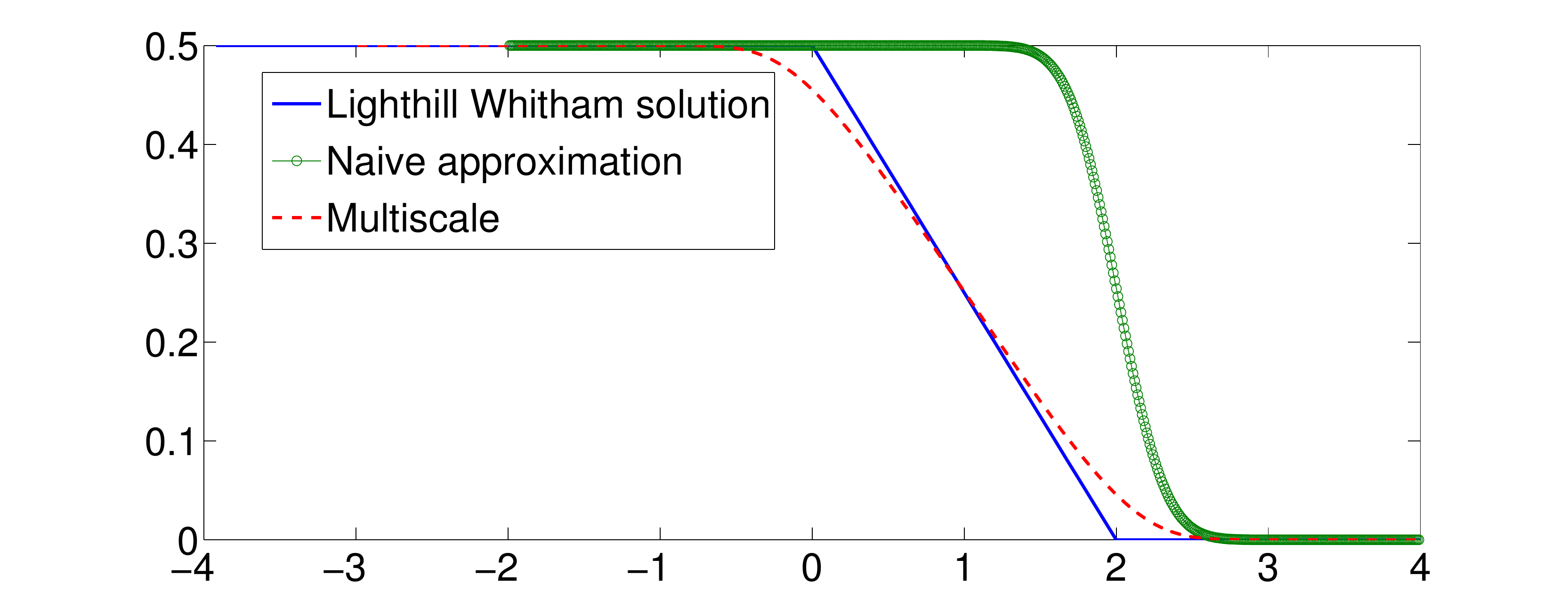}
        }
     \subfloat[$R  = 0.02 $]{
    \includegraphics[keepaspectratio=true, width=.5\textwidth]{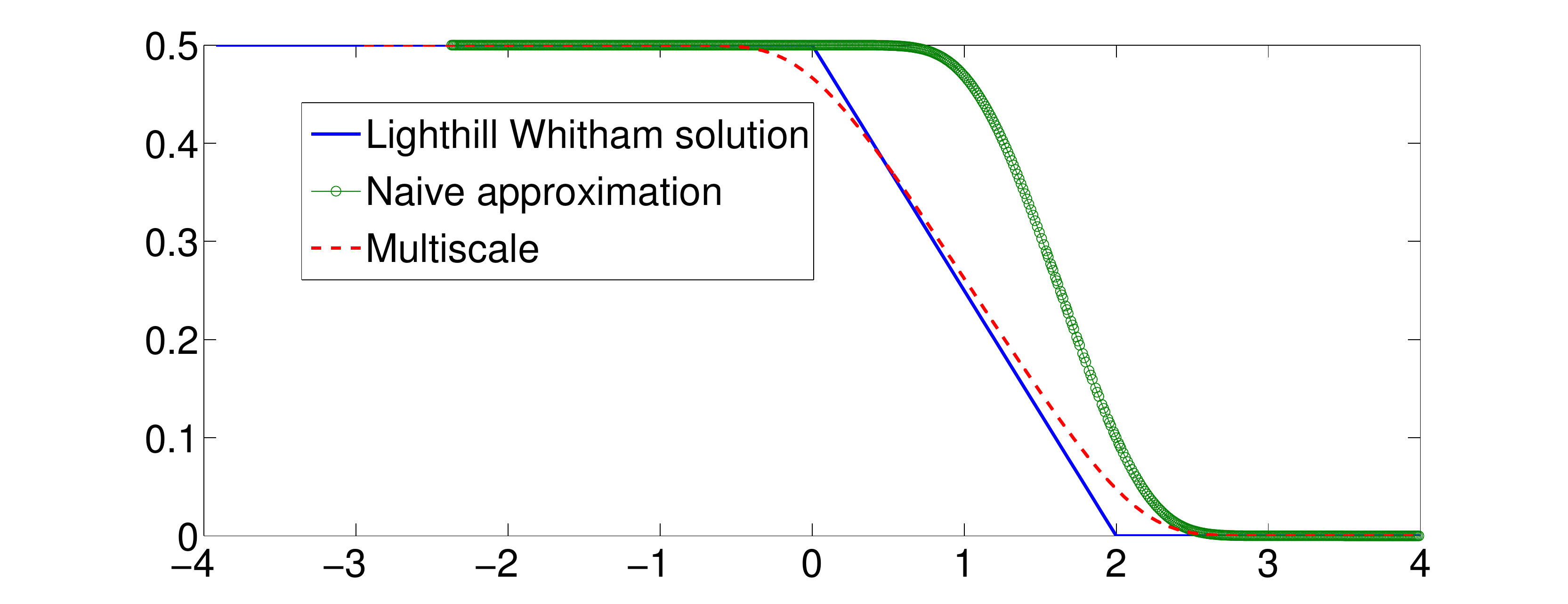}
    }   \\
    \subfloat[ $R = 0.2$]{
        \includegraphics[keepaspectratio=true, width=.5\textwidth]{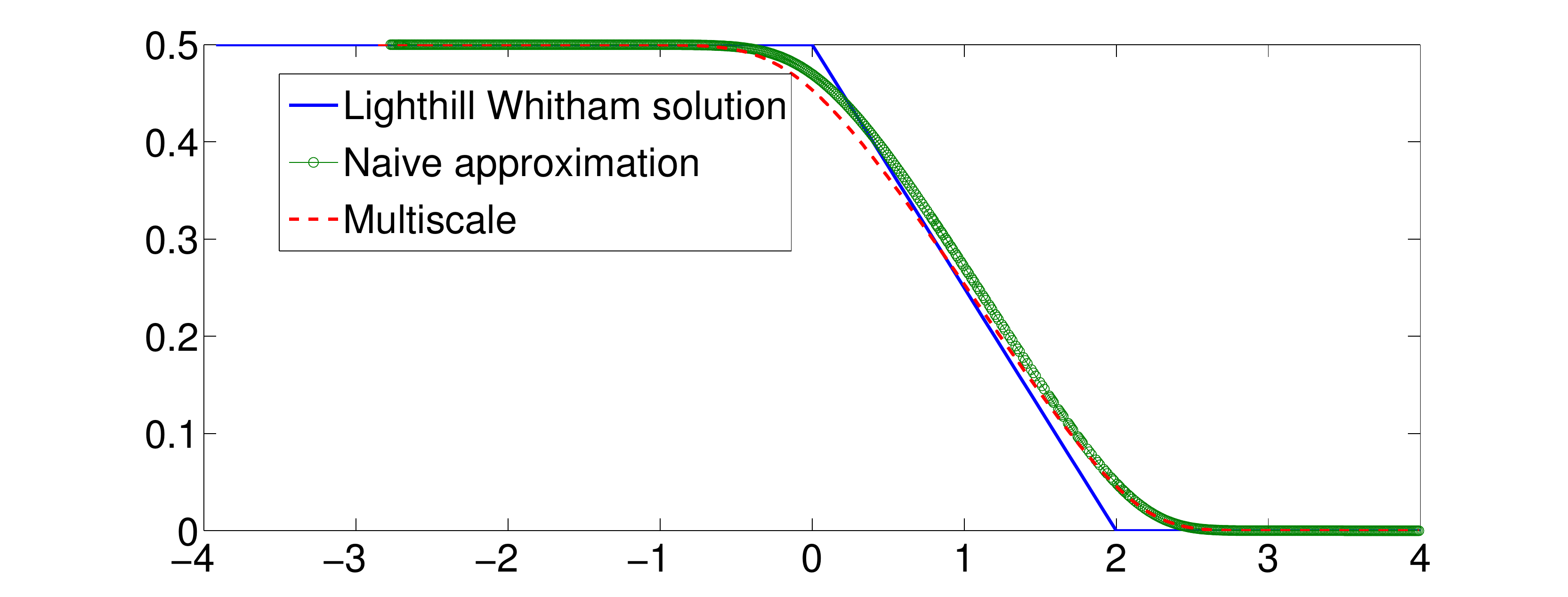}}   
            \subfloat[ $R = 0.4$]{
        \includegraphics[keepaspectratio=true, width=.5\textwidth]{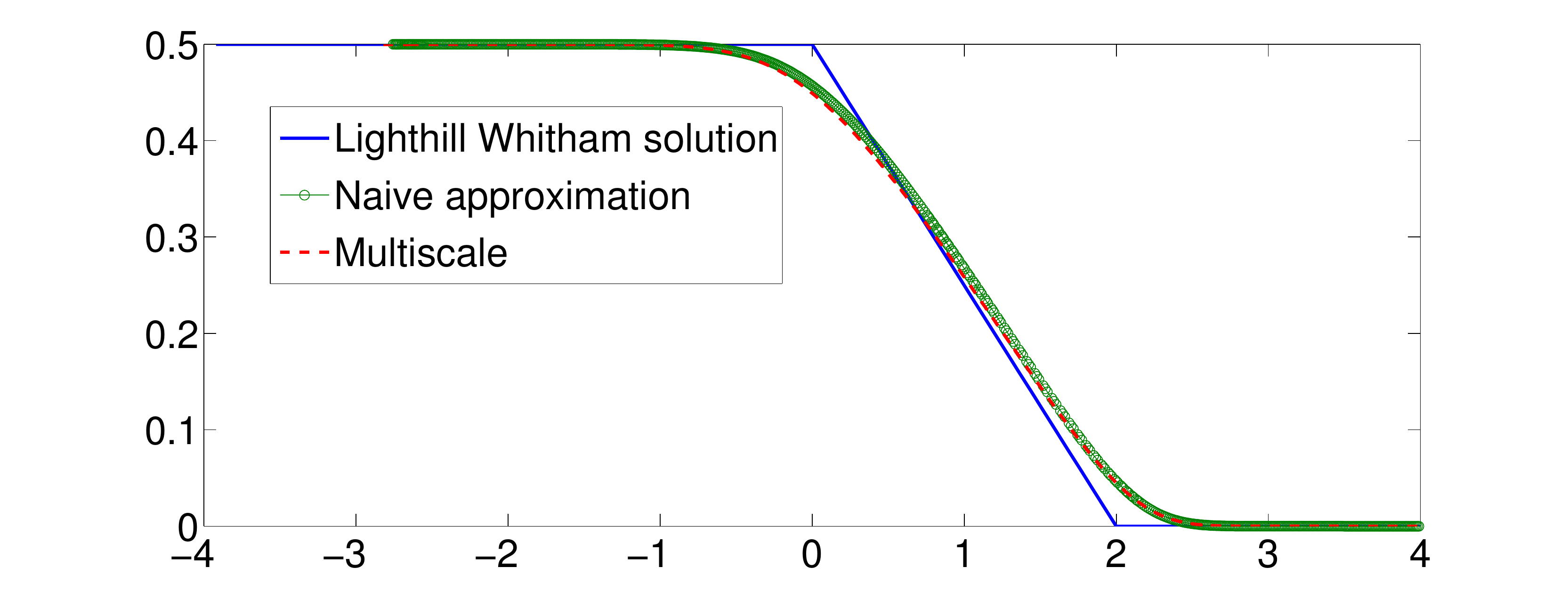}}   
    \captionsetup{margin=0pt}
    \caption{Rarefaction solution  for the  scalar model (\ref{trafficlocal}) and solution of (\ref{trafficnonlocal}) with microscopic and multi-scale interaction approximation at time $t=2$.
    The number of particle is $N=800$ and the interaction radius ranges from $R = 0.002$ to  $R = 0.4$.}
    \label{fig-comparison1}
    \end{figure}

We consider equations (\ref{trafficnonlocal}) and (\ref{trafficlocal}) with different initial conditions and different choices of the interaction potential.
For this example, we consider a fully Lagrangian approach.
  First we consider a symmetric  potential as in Example \ref{quadratic} and initial conditions leading to a  rarefaction wave
  solution. We choose 
 \begin{align}
  \rho_0  (x) =  \begin{cases}
  \frac{1}{2}, \, \,x<0\\
  0, \, \, x>0.
  \end{cases}
  \end{align}
  The solution of the limit equation (\ref{trafficlocal}) is a rarefaction wave 
  \begin{align}
  \rho (x,t)   =  \begin{cases}
  \frac{1}{2}, \, \,x < 0 \\
  \frac{1}{2}\left(1-\frac{x}{t}\right), \, \,0< x < t \\
  0, \, \, x > t .
  \end{cases}
  \end{align}
 We choose and $\delta =0.1$.
  In Figure \ref{fig-comparison1} the rarefaction solution  is compared at time $t=2$ for the  scalar model (\ref{trafficlocal}) and  for (\ref{trafficnonlocal}) with microscopic interaction approximation and multi-scale  approximation. We use a fixed number of  
    particles $N=800$ and an interaction radius ranging from $R = 0.002$ to  $R = 0.4$.
   One observes a good coincidence of microscopic and multi-scale approximation for large $R$ and a stronger deviation, the smaller the value of $R$ is chosen. In this situation the influence of  larger  values of $R$ on the exact solution is a small increase of the smearing of the solution.
   We note that in this situation using a one-sided downwind interaction potential as in Example 2 and $\delta =0$ gives similar results.

\begin{figure}
    \centering
    \captionsetup[subfigure]{margin=0pt} 
     \subfloat[ $R  = 0.002$]{
        \includegraphics[keepaspectratio=true, width=.5\textwidth]{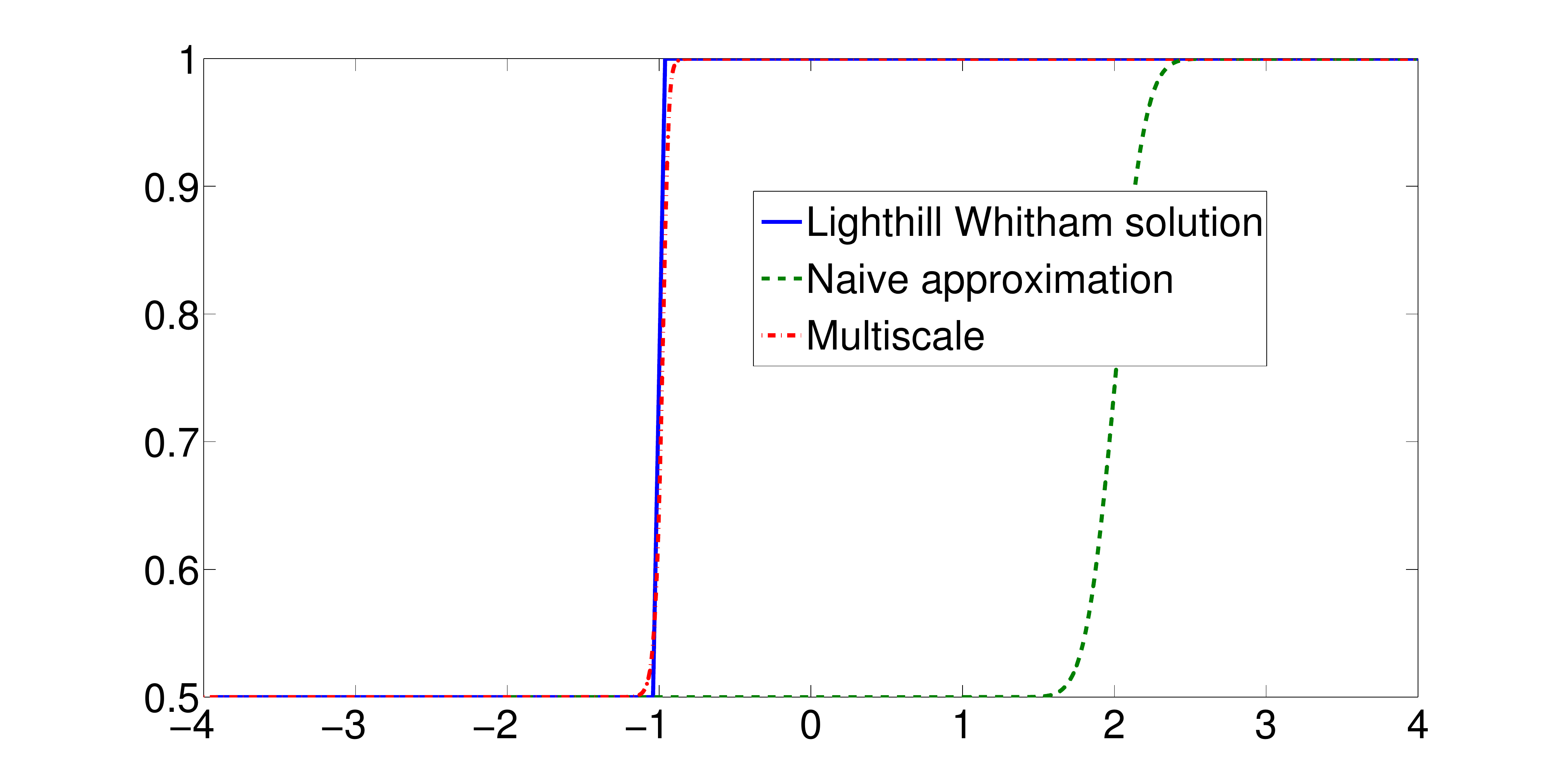}
        }
     \subfloat[$R  = 0.02 $]{
    \includegraphics[keepaspectratio=true, width=.5\textwidth]{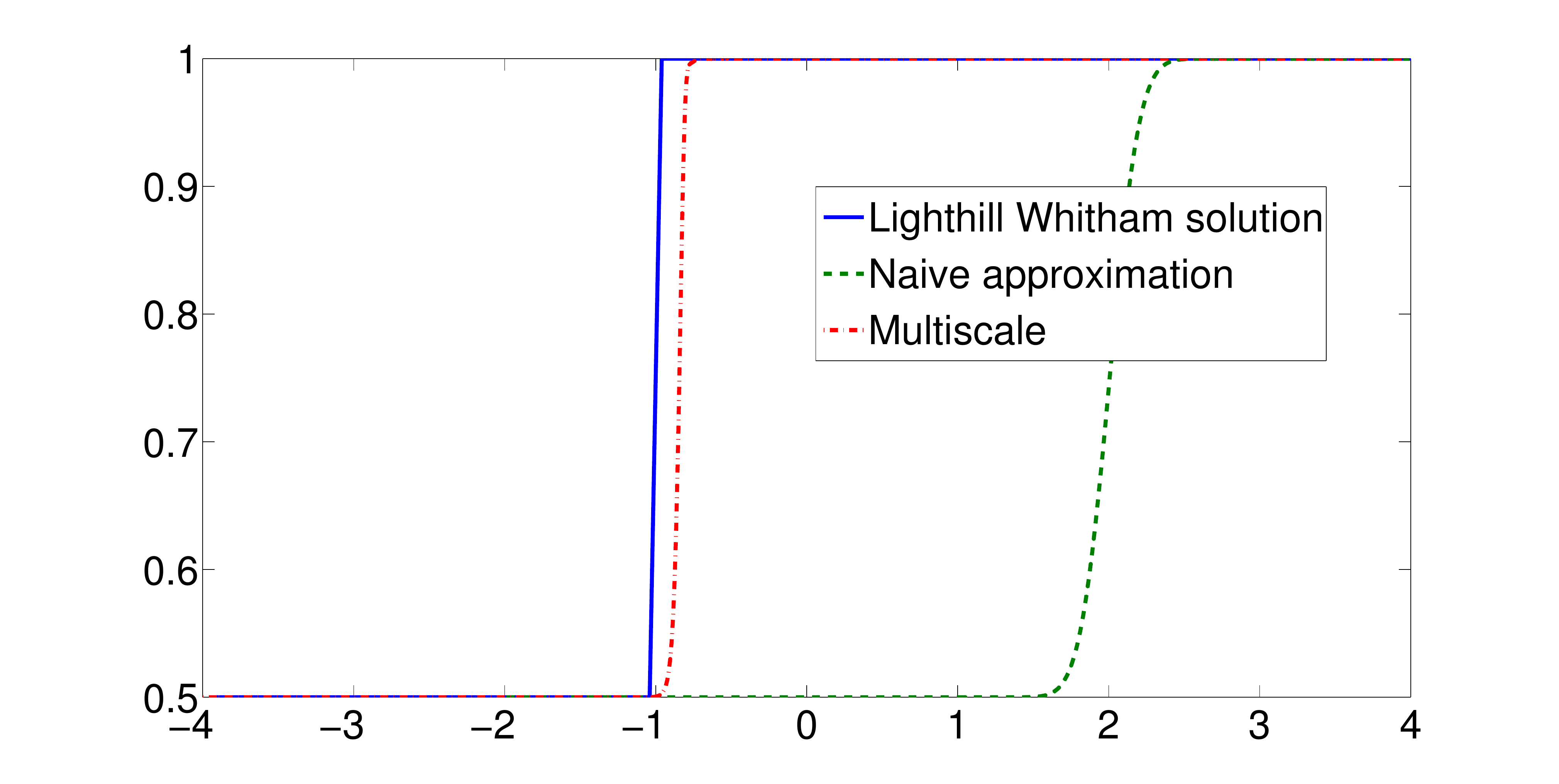}
    }   \\
    \subfloat[ $R = 0.2$]{
        \includegraphics[keepaspectratio=true, width=.5\textwidth]{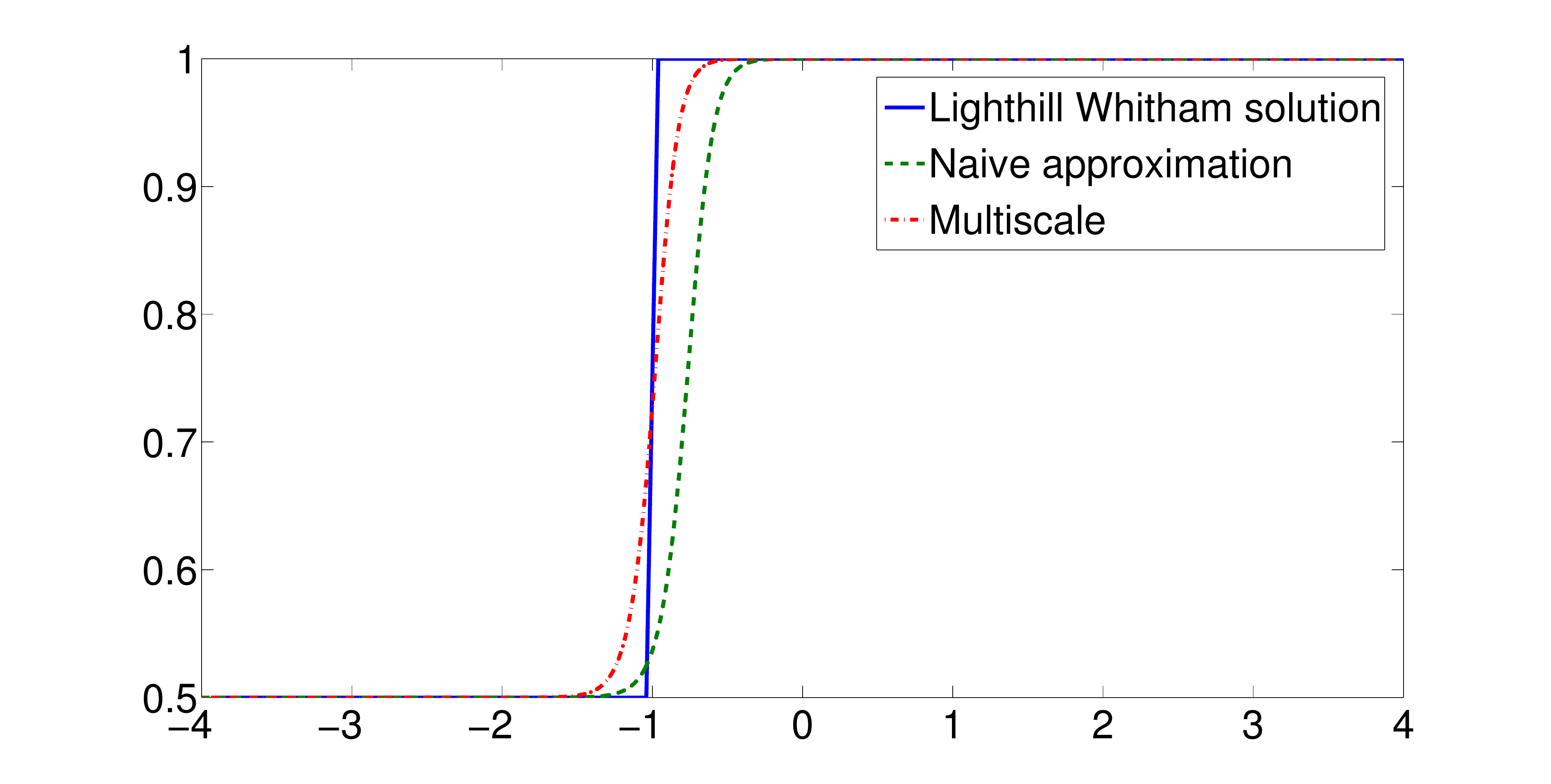}}   
            \subfloat[ $R = 0.8$]{
        \includegraphics[keepaspectratio=true, width=.5\textwidth]{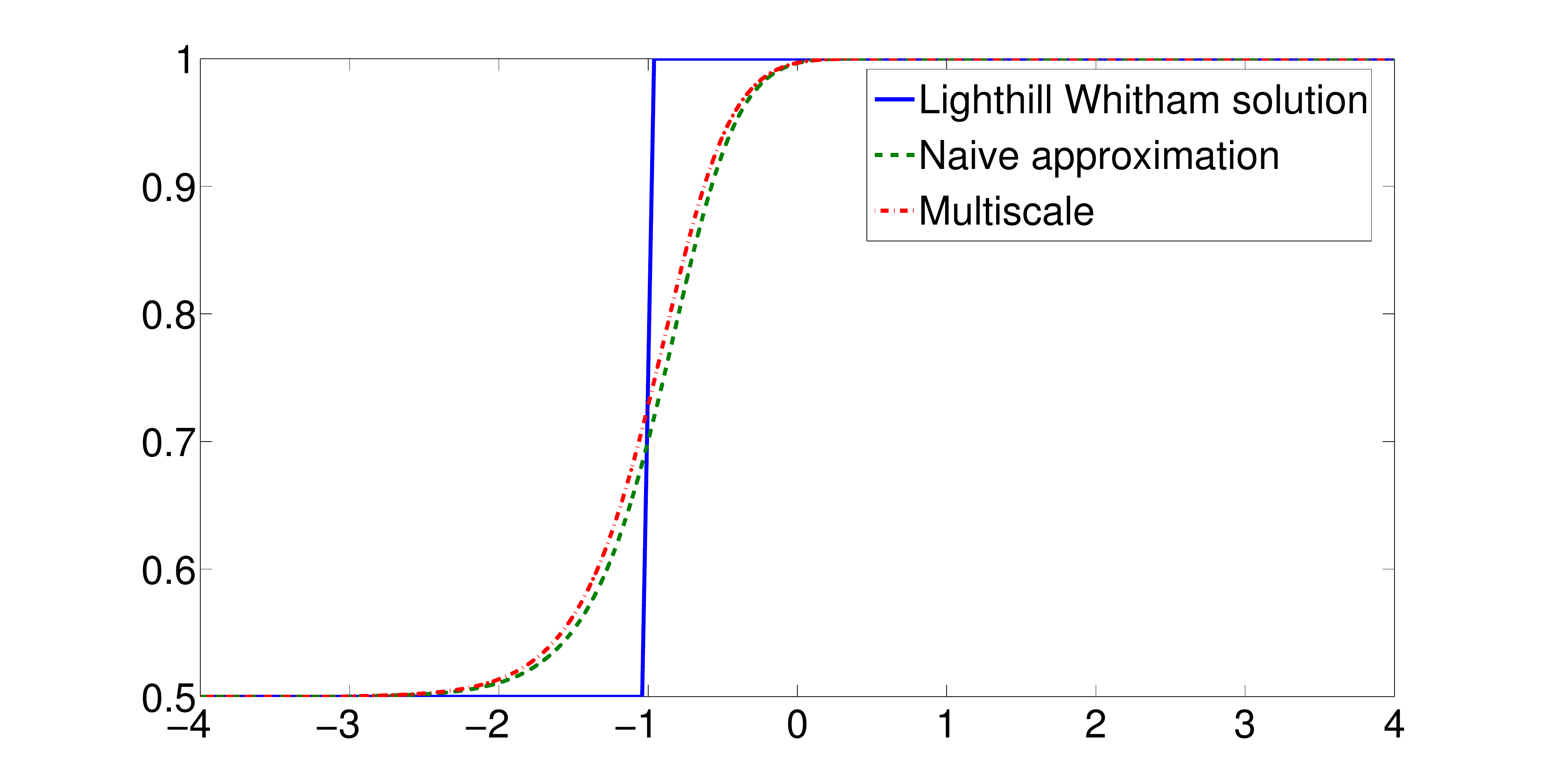}}   
    \captionsetup{margin=0pt}
    \caption{Shock solution for  $N=800$ particles and interaction radius ranging from $R = 0.002$ to  $R = 0.8$  for the  scalar model (\ref{trafficlocal}) and the non-local version  (\ref{trafficnonlocal})  with microscopic and multi-scale interaction approximation  and  downwind potential at time $t=2$ .}
    \label{fig-comparison2}
    \end{figure}

Second  we consider a shock solution of the  Riemann problem for equation (\ref{trafficnonlocal})
and (\ref{trafficlocal}).
 The initial values are now chosen as 
 \begin{align}
 \rho_0  (x) =  \begin{cases}
 \frac{1}{2}, \, \,x<0\\
 1, \, \, x>0.
 \end{cases}
 \end{align}
 The solution of the limit equation (\ref{trafficlocal}) is a shock 
 \begin{align}
 \rho (x,t)   =  \begin{cases}
 \frac{1}{2} , \, \,x <-\frac{t}{2} \\
 1, \, \, x > - \frac{t}{2} .
 \end{cases}
 \end{align}
 We consider first the case of a one-sided downwind potential chosen as in Example 2 and chose $\delta =0$.
 In Figure \ref{fig-comparison2} we compare again   different values of $R$ for a fixed number $N=800$
 of particles. We observe,  that in the well resolved case for larger values of $R$ the microscopic approximation  and the multi-scale solution 
 give a smeared out shock solution deviating from the solution of the limiting equation which is  a moving shock. For the under-resolved case with smaller values of $R$, the microscopic approximation  does not give the correct results, in particular, the speed of the computed wave is wrong and equal to the one of the advection problem obtained after setting $U_R \star \rho$ in (\ref{trafficnonlocal}) equal to $0$.
 This is in  contrast to the  multi-scale method which computes the correct wave speed.

\begin{figure}
          \centering
             \includegraphics[keepaspectratio=true, width=.5\textwidth]{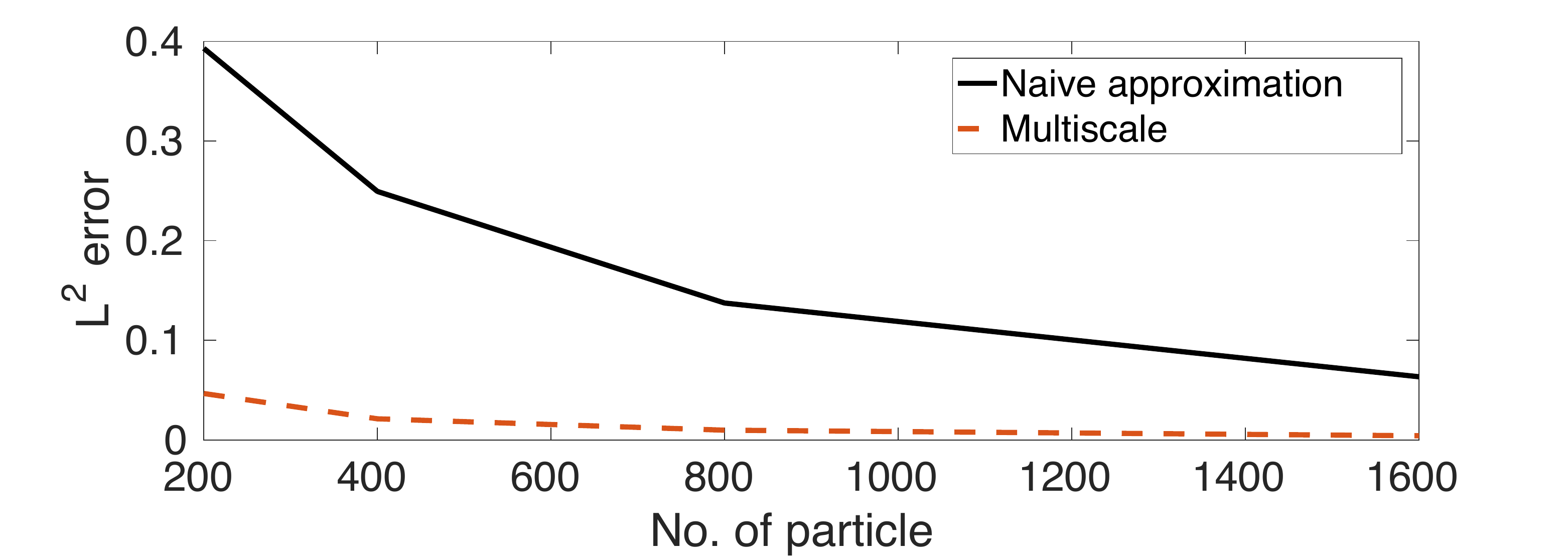}
             \caption{$\mathcal{L}^2$-error plot at time $t=2$  computed from    (\ref{trafficnonlocal})  
                                   with $R= 0.2$ using the microscopic and the multi-scale  approximation and a  downwind potential.}
             \label{error2}
             \end{figure}

For the investigation of  the convergence of the method 
  we consider  a situation  with $R=0.2$. The scalar problem is solved with microscopic and multi-scale approximation of the interaction term.
  We compute the error of the different  methods for different numbers of grid-particles comparing the numerical results to the numerical solution of (\ref{trafficnonlocal}) determined with a fine grid with 
  $N=6400$. For  $6400$  particles the difference between multi-scale and naive method is again of the order $10^{-2}$.  In Figure \ref{error2} and Table \ref{table2} the $\mathcal{L}^2$-errors  versus the number of particles are plotted.       
   One observes again the deterioration of the   method using the microscopic approximation of the interaction term for smaller numbers of particles.
     The multi-scale method is able to treat all ranges with a reasonable accuracy for  the first order method.

 \begin{table}
 \begin{center}
 \begin{tabular}{|r|r|r|r|}
  \hline
 $\# $ particles  & naive  &    multi-scale    & CPU time   \\ 
 &error&error& in seconds  \\\hline 
   $200$ & $0.40$    & $0.05$   & 8  \\
 $400$ &$0.26$   & $0.02$     & 17  \\ 
 $800$  &$0.15$  & $0.01$    &  36\\
 $1600$ &$0.08$   & $0.002$    & 77  \\
 \hline
 \end{tabular}
 \caption{Convergence study for nonlocal Lighthill-Whitham equations with downwind interaction potential and 
 $R=0.2$.}
 \label{table2}
 \end{center}
 \end{table}

Similiar to  the  example in the last subsection comparable errors are obtained using $N=200$ particles and a computation time of $8s$ for the multi-scale method and $N=1600$ particles and a computation time of $77s$ for the microscopic method. This yields again an order of magnitude gain in computation time.

\begin{figure}
    \centering
    \captionsetup[subfigure]{margin=0pt} 
     \subfloat[ $R  = 0.002$]{
        \includegraphics[keepaspectratio=true, width=.5\textwidth]{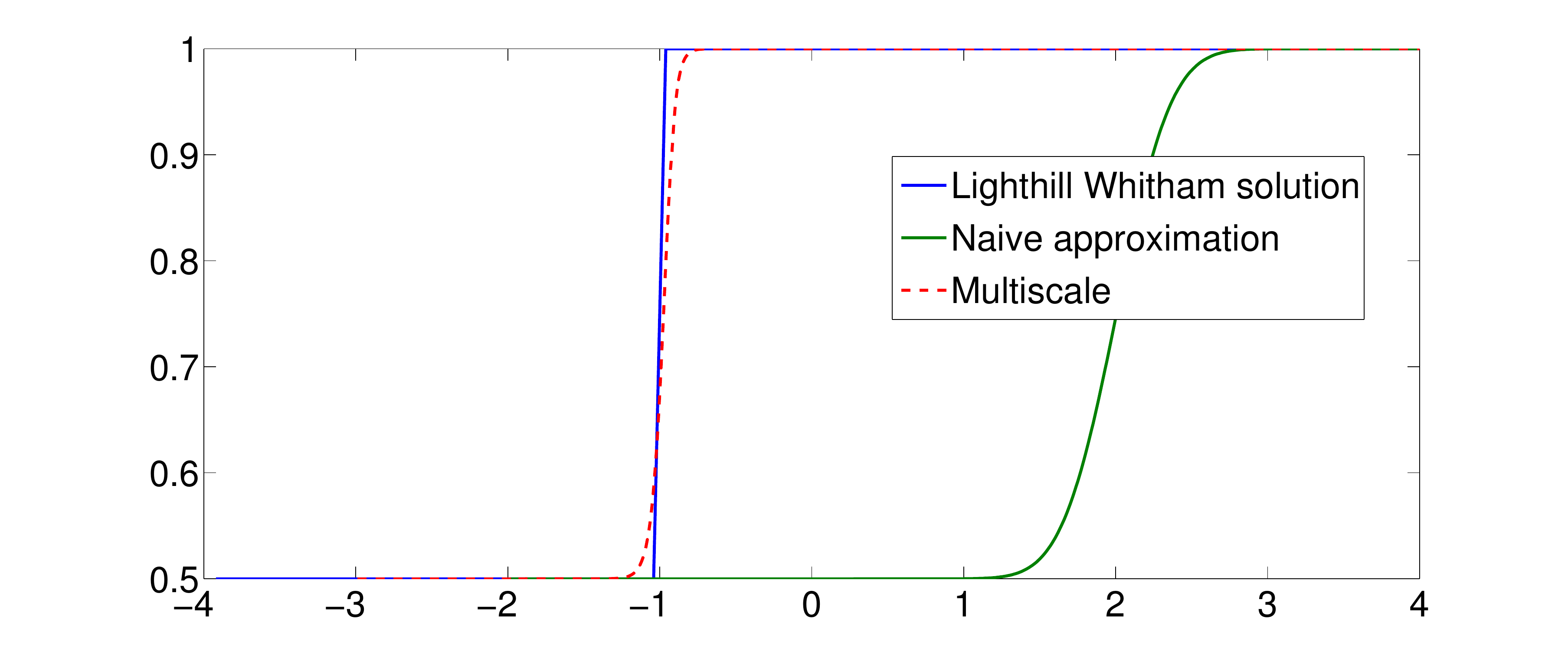}
        }
     \subfloat[$R  = 0.02 $]{
    \includegraphics[keepaspectratio=true, width=.5\textwidth]{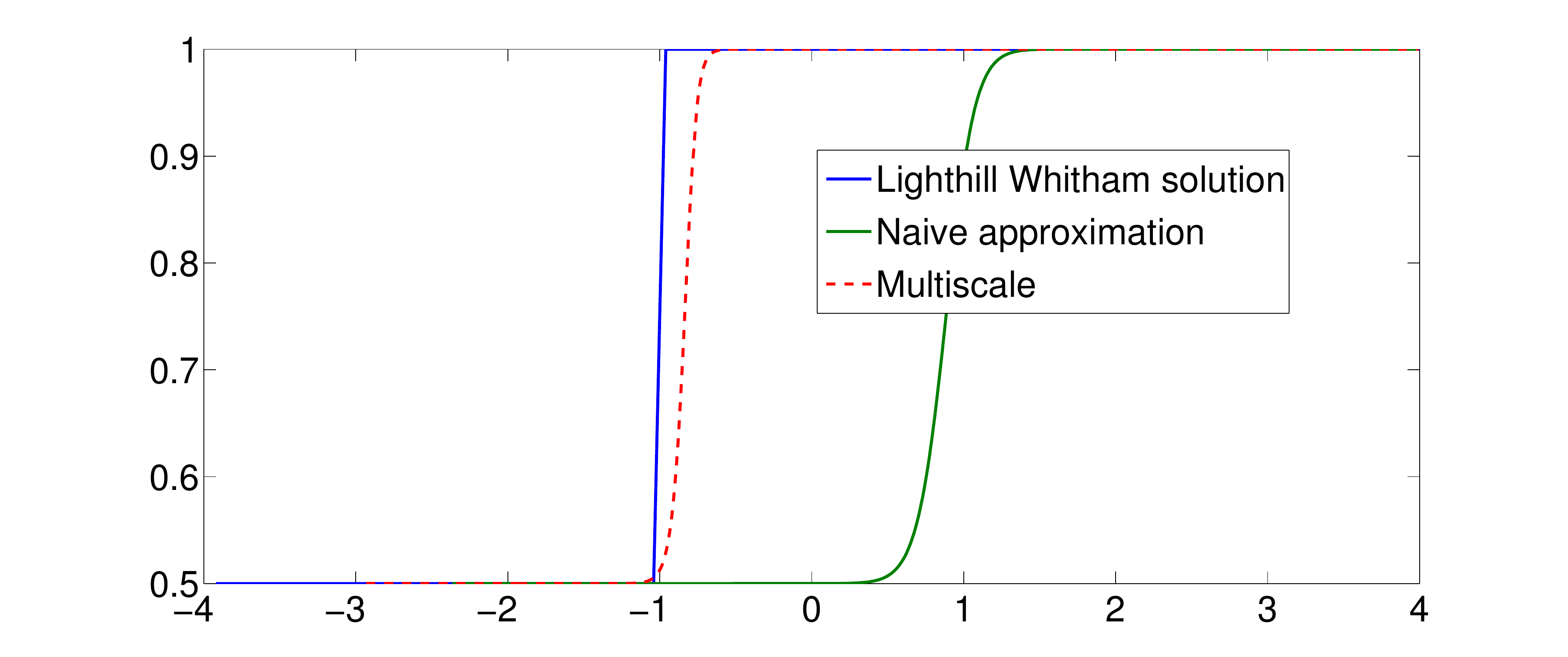}
    }   \\
    \subfloat[ $R = 0.2$]{
        \includegraphics[keepaspectratio=true, width=.5\textwidth]{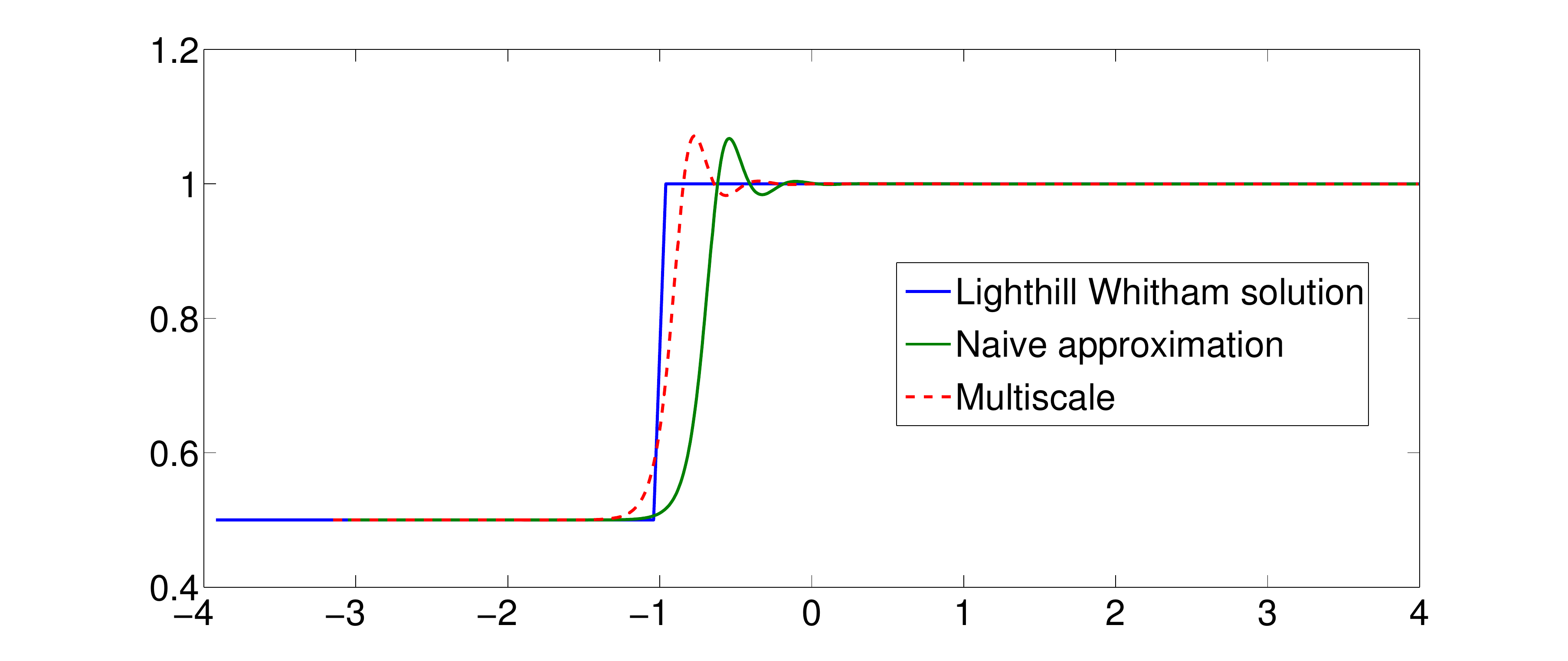}}   
            \subfloat[ $R = 0.8$]{
        \includegraphics[keepaspectratio=true, width=.5\textwidth]{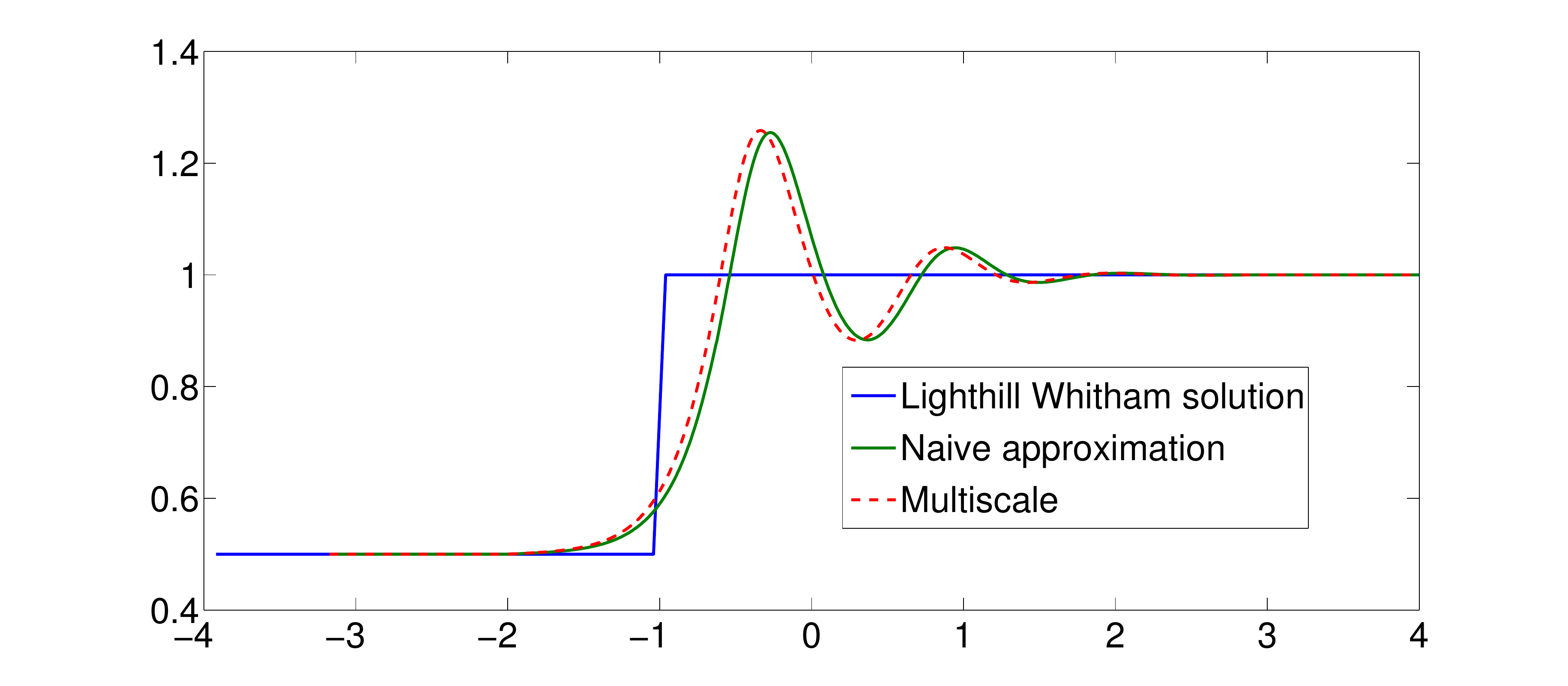}}   
    \captionsetup{margin=0pt}
    \caption{Shock solution for  $N=800$ particles and $R = 0.002$ to  $R = 0.8$ for the  scalar model (\ref{trafficlocal}) and the non-local version  (\ref{trafficnonlocal})  with microscopic interaction approximation and multi-scale  approximation at time  $t=2$ for the symmetric potential.}
    \label{fig-comparison2symm}
    \end{figure}

Finally we investigate the above example with a symmetric interaction potential and a regularization $\delta =0.02$. 
 In Figure \ref{fig-comparison2symm} we compare again   different values of $R$ for a fixed number $N=800$
of particles. We observe,  as before that in the well resolved case for larger values of $R$ the microscopic approximation  and the multi-scale solution 
give an oscillating  solution deviating from the  moving shock solution of the limit equation. For the under-resolved case with smaller values of $R$, the microscopic approximation  does not give the correct results as in the case of the  downwind potential, whereas
the  multi-scale method  computes the correct solution.

\begin{figure}
          \centering
             \includegraphics[keepaspectratio=true, width=.5\textwidth]{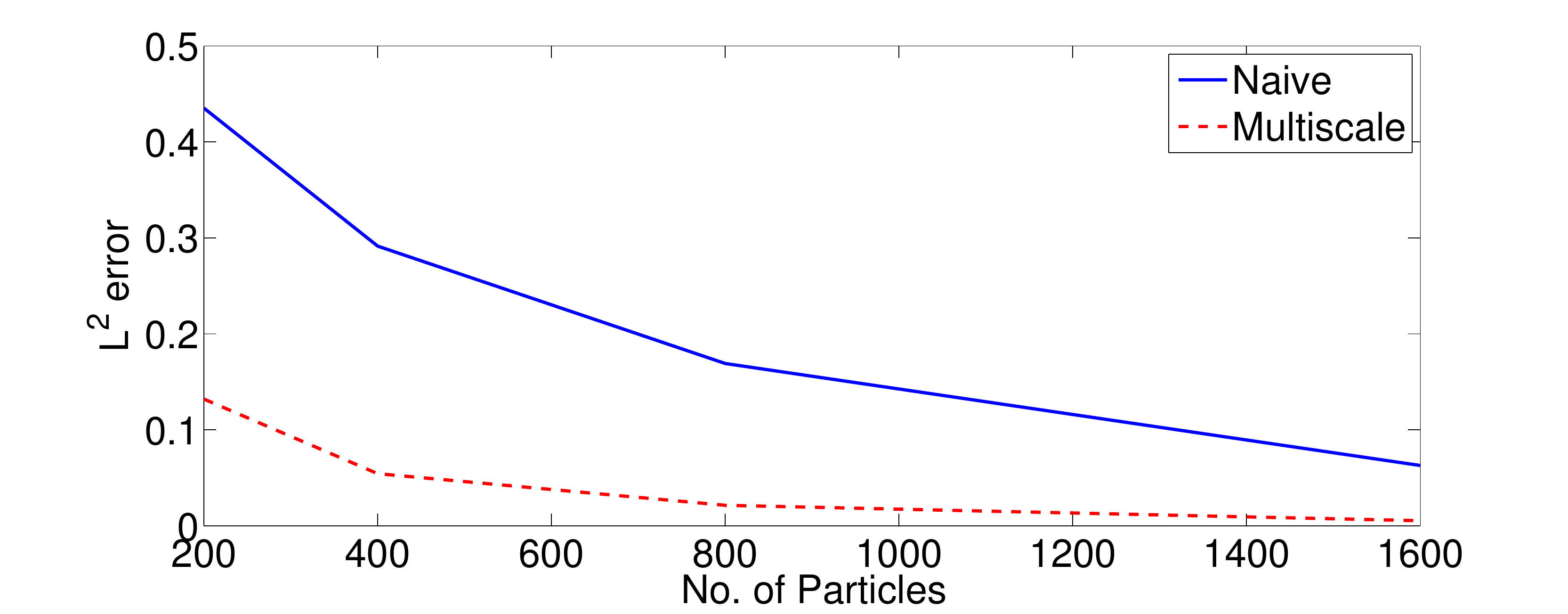}
             \caption{$\mathcal{L}^2$-error plot for the solutions at time $t=2$  computed from    (\ref{trafficnonlocal})  
                                   with $R= 0.1$ and $\delta =0.02$ using the microscopic and the multi-scale  approximation for symmetric interaction potential.}
             \label{error2a}
             \end{figure}

For the investigation of  the convergence of the method 
we consider as before  $R=0.2$.
We compute the error of the particle method for different numbers of grid-particles comparing the numerical results to the numerical solution of (\ref{trafficnonlocal}) determined with a fine grid with 
$N=6400$. The scalar problem is solved with microscopic and multi-scale approximation of the interaction term. In Figure \ref{error2a} and Table \ref{table2a} the $\mathcal{L}^2$-errors  versus the number of particles are plotted.       
One observes a similiar behaviour as in the non-symmetric potential case.
We  remark that the method developed here is still working in the present case, since we have regularized the equations with $\delta >0$ obtaining a  dissipative limit. For $\delta $ much smaller,
the limit would be dominated by dispersion and the method is not supposed to work more efficienly than the naive, microscopic approximation.

\begin{table}
 \begin{center}
 \begin{tabular}{|r|r|r|r|}
  \hline
 $\# $ Particles  & naive  &    multi-scale    & CPU time   \\ 
 &error&error&in seconds  \\\hline 
   $200$ & $0.45$    & $0.13$   & $9$   \\
 $400$ &$0.31$   & $0.05$     &   $19$ \\ 
 $800$  &$0.21$  & $0.02$    &  $50$ \\
 $1600$ &$0.11$   & $0.01$    &  $142$   \\ \hline
 \end{tabular}
 \caption{Convergence study for nonlocal Lighthill-Whitham equations with symmetric interaction potential and 
 $R=0.2$ and $\delta =0.02$.}
 \label{table2a}
 \end{center}
 \end{table}

  \subsection{Comparison of numerical algorithms  in 2D}

 In this section  we solve the hydrodynamic system (\ref{pednonlocalhydro}) using the  multi-scale  algorithm and the  
  algorithm with naive, microscopic evaluation of the integral  in 2D.  
 We compare the results  with each other and with the  local hydrodynamic 
 problem  (\ref{pedlocalhydro}). We consider a fully Lagrangian approach.

  \subsubsection{Hydrodynamic system for test case}
 
\begin{figure}
    \centering
    \captionsetup[subfigure]{margin=0pt} 
     \subfloat[]{
        \includegraphics[keepaspectratio=true, width=.5\textwidth]{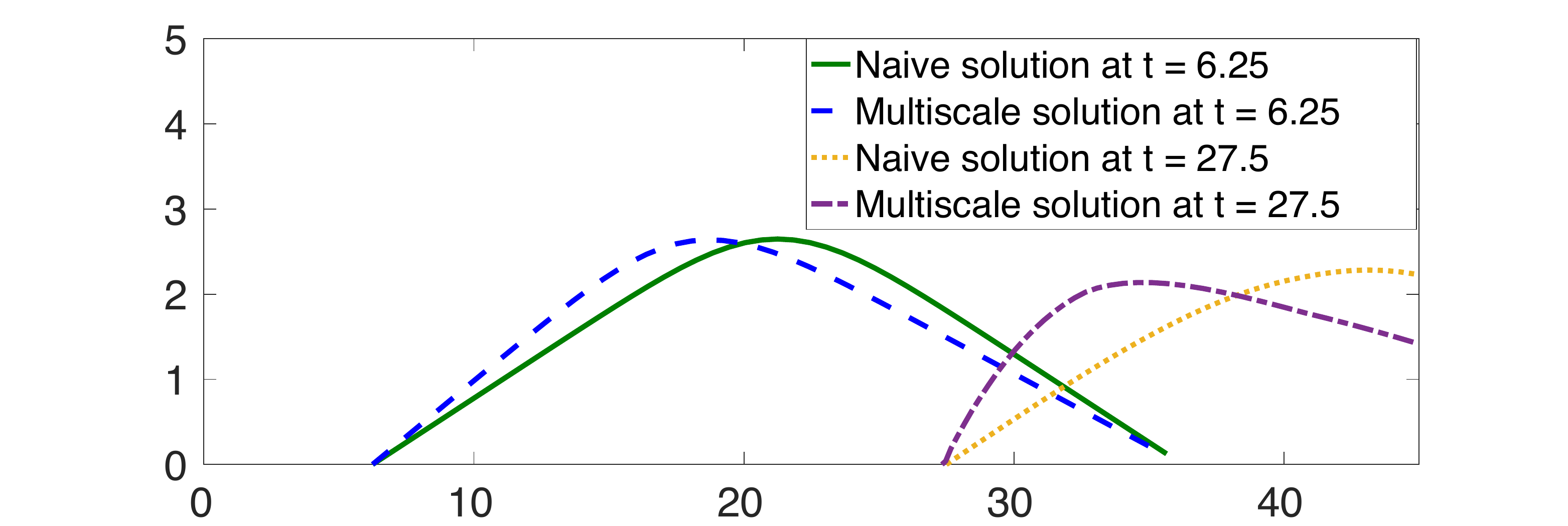}
        }
     \subfloat[]{
    \includegraphics[keepaspectratio=true, width=.5\textwidth]{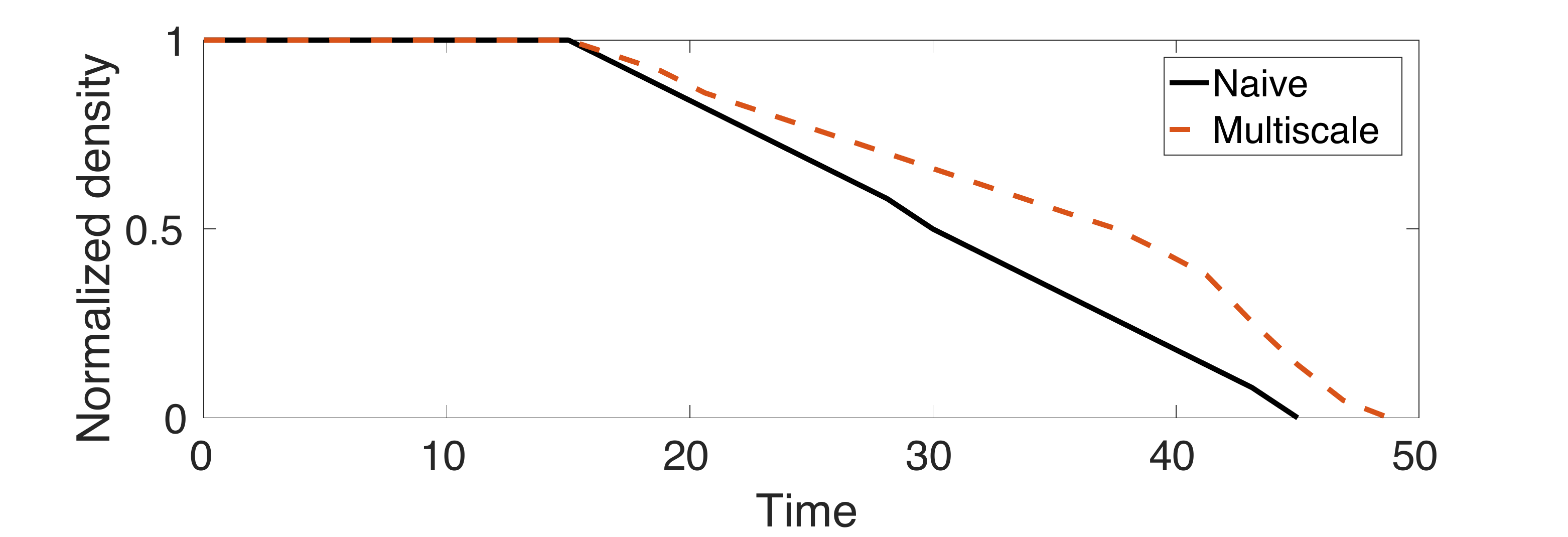}
    }  
    \captionsetup{margin=0pt}
    \caption{Solution for  initial spacing $\Delta x = 0.6$ or $N=4150$ particles and $R = 0.2$. (a) Solutions  for the  hydrodynamic model  (\ref{pednonlocalhydro})   with microscopic interaction approximation and multi-scale  approximation at time  $t=6.25$ and $27.5$. (b)Time development of the normalized total mass in the computational domain determined from the hydrodynamic pedestrian model  (\ref{pednonlocalhydro}) using  microscopic  and  multi-scale  approximation for fixed interaction radius $R= 0.2$.}
    \label{test}
    \end{figure}

 We consider a two dimensional domain $[0,45]\times[0, 50] $.  We initially generate particles in $[0, 30]\times[0, 50]$. 
 We simplify  the hydrodynamic system (\ref{pednonlocalhydro}) considering a fixed vector 
 $
 \hat e (x) = (1,0)^T.
 $
 The maximal velocity is $u_{max}=1$ and  the maximal density is chosen as $\rho_{max}=1$. 
 A symmetric interaction potential and $\delta = 0.1$ are chosen. Moreover, $\alpha =1, \gamma =1$.  We use  $N=4150$ particles  and the interaction radius $R$ equal to $0.2$.  
The  initial condition  is  given by 
\begin{align}
 \rho_0  (x) =  \begin{cases}
 \frac{1}{5}x, \, \,0<x \le 15\\
  3-\frac{1}{5} (x-15), \, \,15<x \le 30\\
 0, \, \,30 <  x \le 45.
 \end{cases}
 \end{align}
In Figure \ref{test} a.) the solutions  of the nonlocal equations (\ref{pednonlocalhydro}) with microscopic and multi-scale interaction aproximation for the interaction radius $R = 0.2$ are plotted. 
Moreover, the time development of the normalized total mass in the  domain $x \in [0,100]$ is plotted  for microscopic and multi-scale interaction approximation
in Figure \ref{test} b.).

\begin{figure}
     \centering
        \includegraphics[keepaspectratio=true,   width= .65\textwidth]{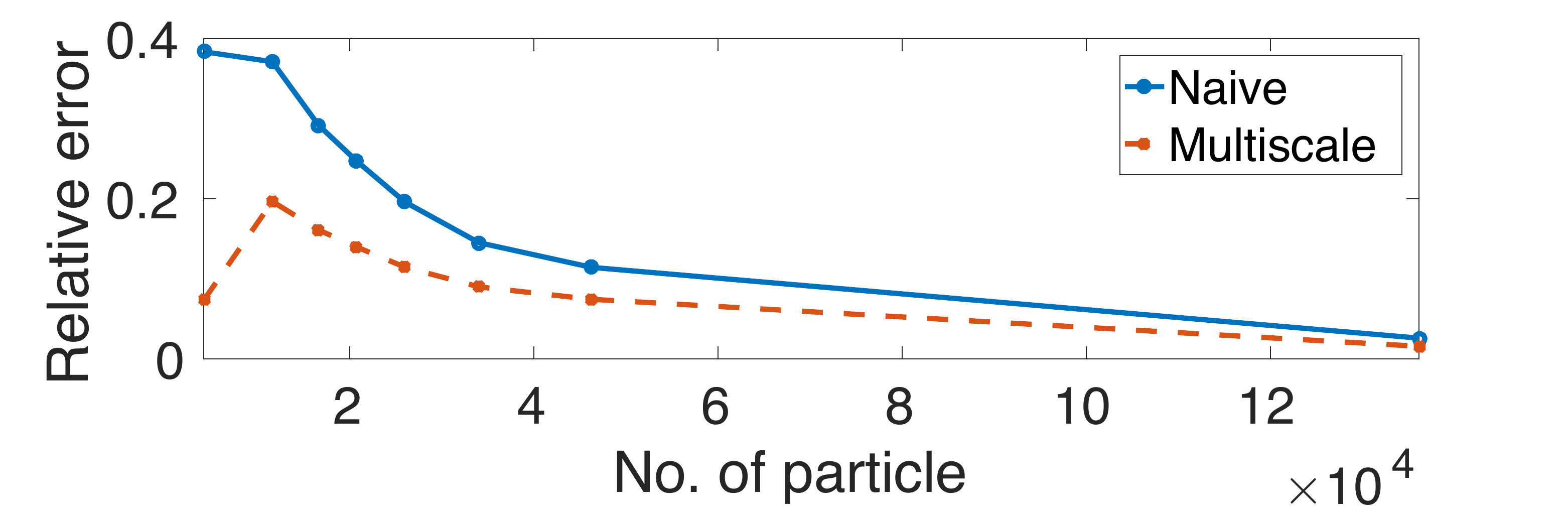}  
        \caption{Error plot for microscopic and multiscale simulations for  hydrodynamic $2D$ equations.}
              \label{fig11}
        \end{figure}

In Table \ref{table4} 
     we compute the error and the CPU times of naive and multi-scale  method for different numbers of grid-particles, see also Figure \ref{fig11}.  The reference solution is computed by  using a spacing of $\Delta x = 0.11$ and
     approximately $136000$ particles. For this very fine resolved case the difference of naive and multi-scale solution is 
     approximately equal to $10^{-2}    $.
       One observes in this physical situation   a smaller but still relevant gain in computation time.
       
 \begin{table}
 \begin{center}
 \begin{tabular}{|r|r|r|r|r|}
  \hline
 initial   &$\# $ particles  & naive &multi-scale & CPU time   \\ 
  spacing &$ $ &error & error  &   \\ \hline 
   $0.6$ & $4150$ &$ 0.38$ & $ 0.07 $   & $$ 15 min   \\
  $0.36$ & $11600$ &$0.37 $ & $0.19 $   & $40$ min   \\
 $0.3$ &$16600$ &$ 0.29$ &$0.16 $   & $71$ min      \\ 
 $0.24$ &$26000$  &$0.20 $ & $0.11$  & $90$ min    \\
 $0.21$ & $34000$ & $0.14$ & $0.09$ & $118$ min \\
 $0.18$ & $46000$ & $0.11$  & $0.07$ & $154$ min \\
 $0.1$ & $136000$ &$ 0.02$ & $0.01$ & $603$ min\\   
 \hline
 \end{tabular}
 \caption{Comparison of CPU times between microscopic and multiscale simulations of the hydrodynamic 2D equations. The error analysis is performed at time 18.75 sec.}
 \label{table4}
 \end{center}
 \end{table}

\subsubsection{Hydrodynamic systems for pedestrian flow (Example 5)}
 In our final example we investigate the hydrodynamic systems  (\ref{pednonlocalhydro})  and (\ref{pedlocalhydro}). Compared to the last subsection we consider a more complicated configuration
defined in Ref. \cite{LWZSL09}, compare also (\cite{EGKT14}), and a full coupling to the eikonal equation.

\begin{figure}
          \centering
             \includegraphics[keepaspectratio=true, width=.5\textwidth]{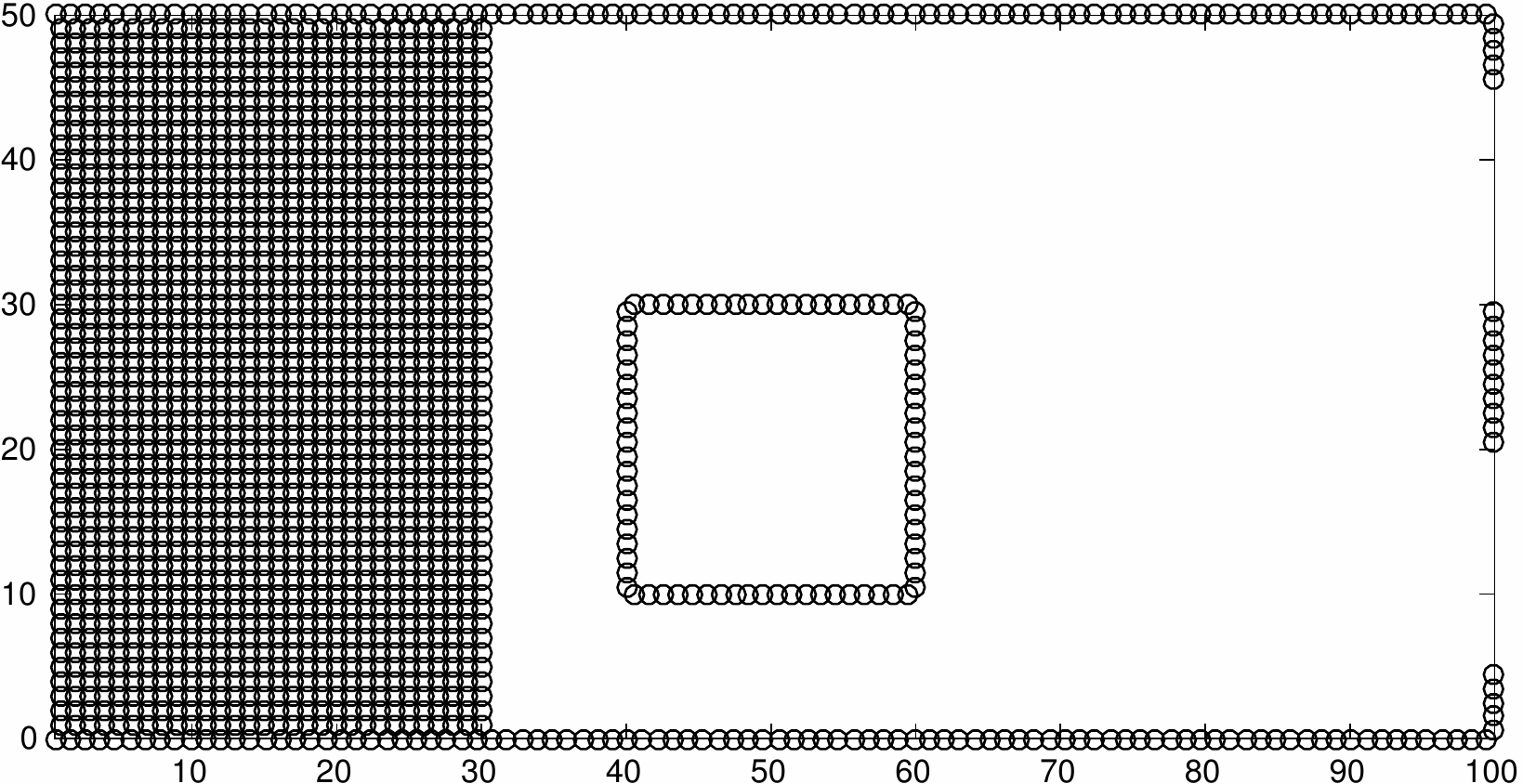}    \hfill
              \caption{ Initial configuration}
             \label{initial_ped}
             \end{figure}

Pedestrians are   initialized on the left of the domain and evacuated towards the exits on the right as shown in 
Figure \ref{initial_ped}. As initial value we choose a constant  value of $\rho = 1  $ in the region
$[0,30] \times [0,50]$.
  In the center  of the computational domain an obstacle is located. For the eikonal equation we use $\phi=0$ on the two exits and 
$\phi = \infty$ on all walls as boundary conditions. 

 We choose the maximum velocity  $u_{max} = 2$ and the maximum density $\rho_{max} = 10  $.   In this case, we  vary the initial average distance between grid points from $0.15$ to $1$, i.e. the number of grid particles varies between $1440$ and $61800$.
 Moreover,  we choose the following parameter:  $\gamma=  500 , \alpha = 1000$.

\begin{figure}
   	\centering
   	  	\captionsetup[subfigure]{margin=5pt} 
   	  \subfloat[localized ]{
   	  	\includegraphics[keepaspectratio=true, width=.5\textwidth]{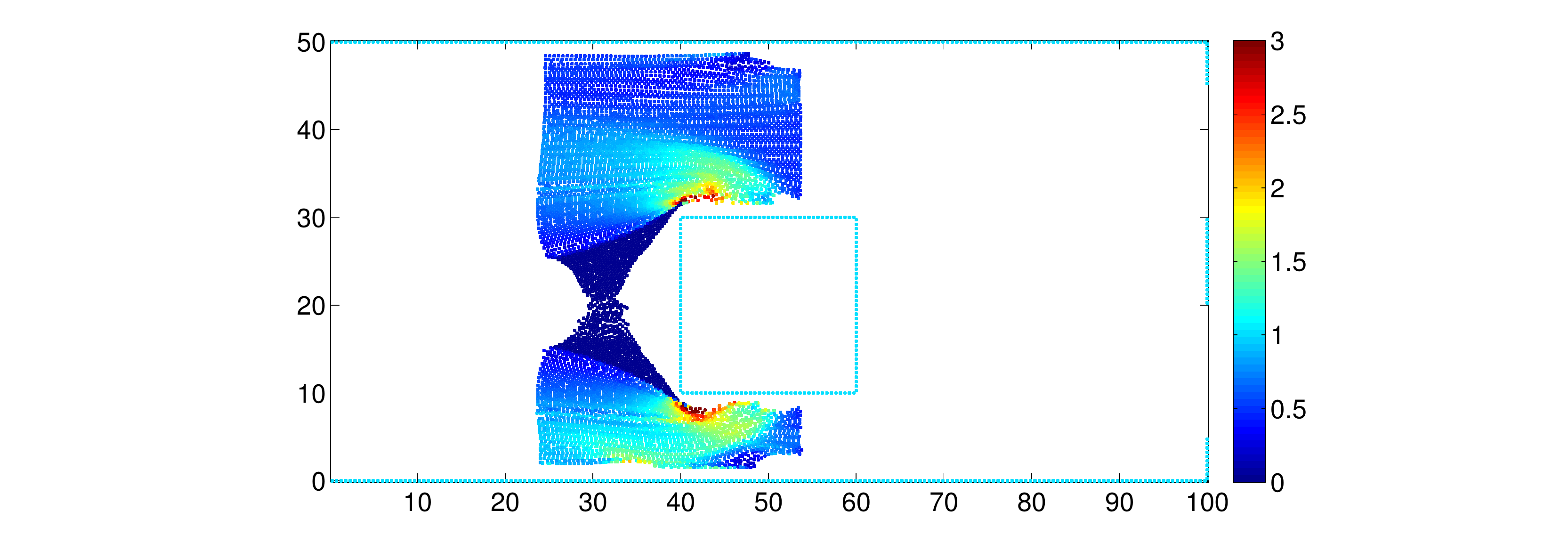}
   	  } \subfloat[ microscopic]{
   	  	\includegraphics[keepaspectratio=true, width=.5\textwidth]{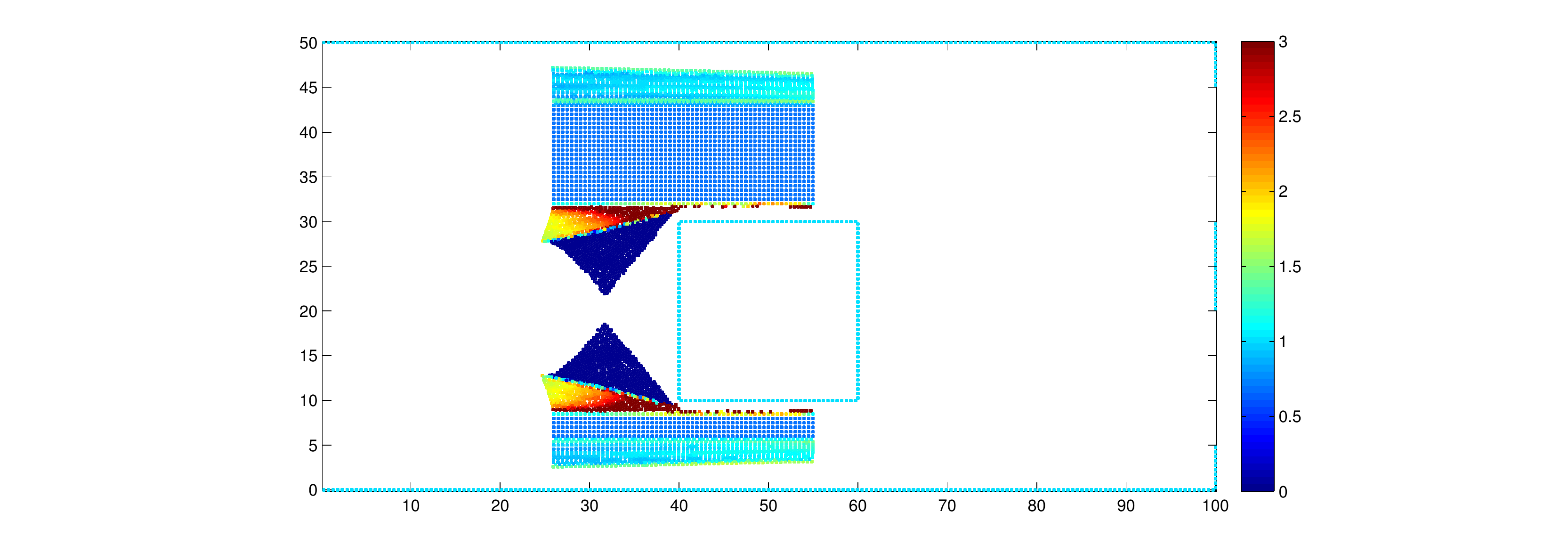}
   	  } 
     \captionsetup[subfigure]{margin=5pt} 
     \subfloat[multi-scale]{
     	\includegraphics[keepaspectratio=true, width=.5\textwidth]{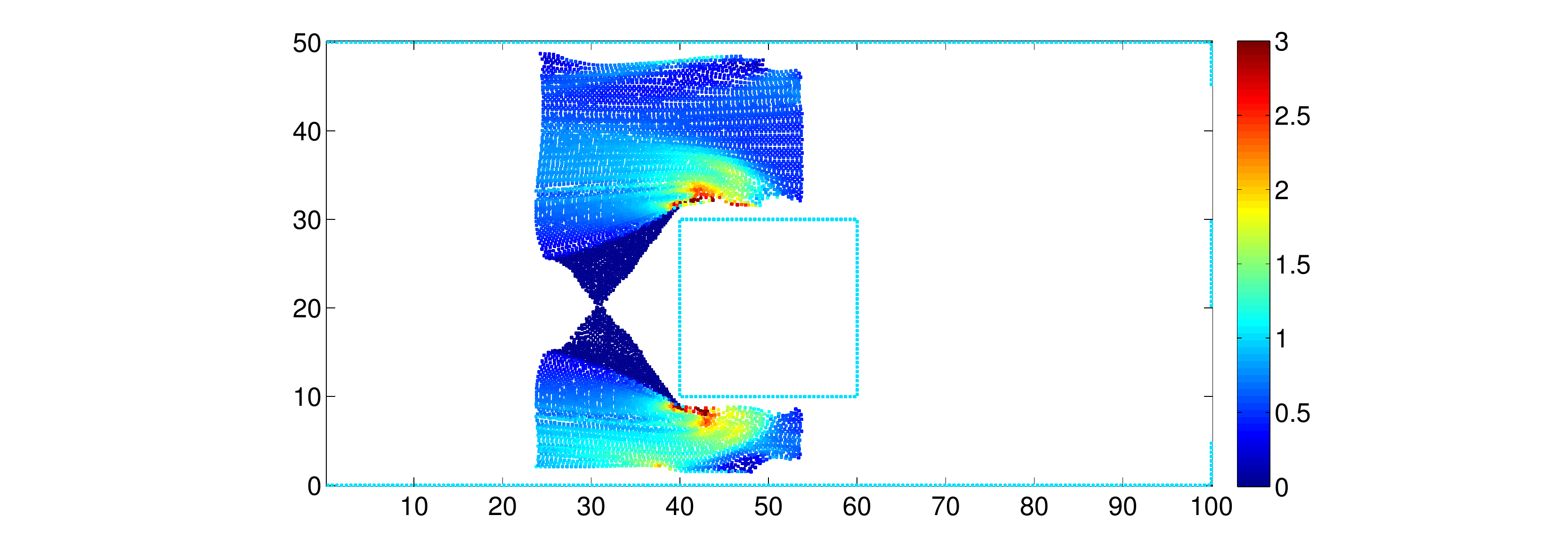}
     } 
              \caption{Density plot determined from local limit equation (\ref{pedlocalhydro})  and nonlocal equations (\ref{pednonlocalhydro}) with microscopic and multi-scale approximation at time $t = 12.5$ for $\Delta x = 0.5$,  $R = 0.2$.}
             \label{compare_ped1}
             \end{figure}

\begin{figure}
          \centering
          	\captionsetup[subfigure]{margin=5pt} 
          	\subfloat[$\Delta x =0.4$]{
          		\includegraphics[keepaspectratio=true, width=.5\textwidth]{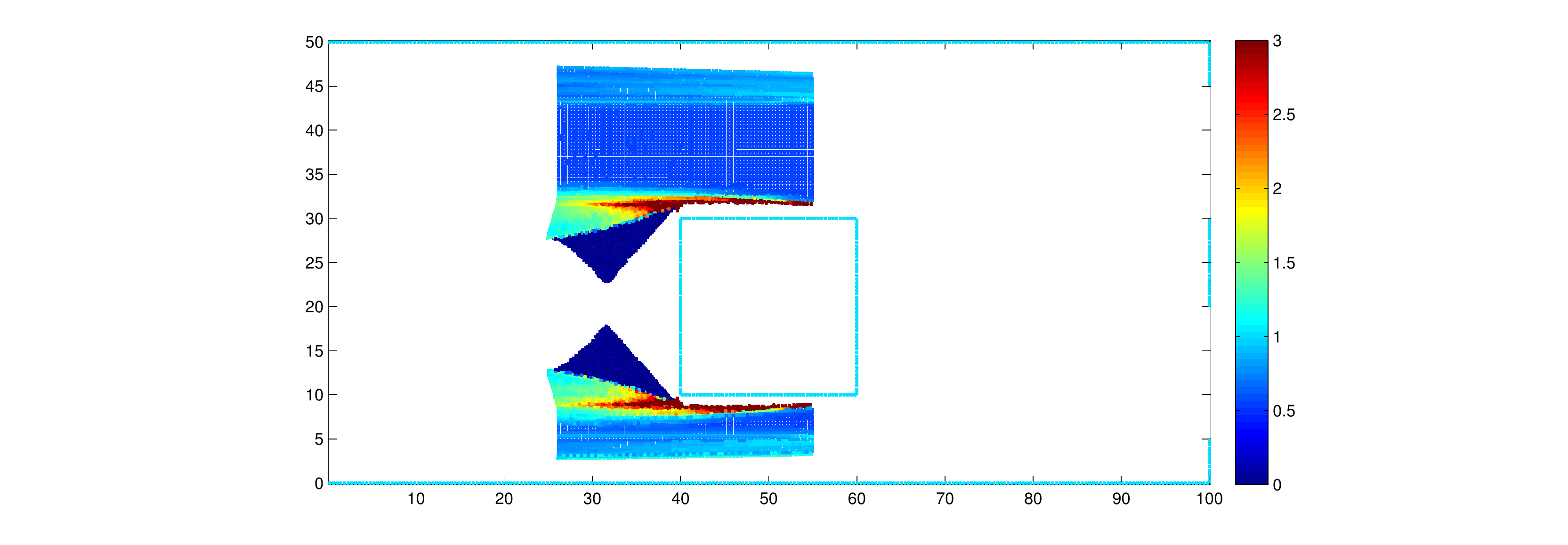}
          	} \subfloat[ $\Delta x = 0.3$]{
          		\includegraphics[keepaspectratio=true, width=.5\textwidth]{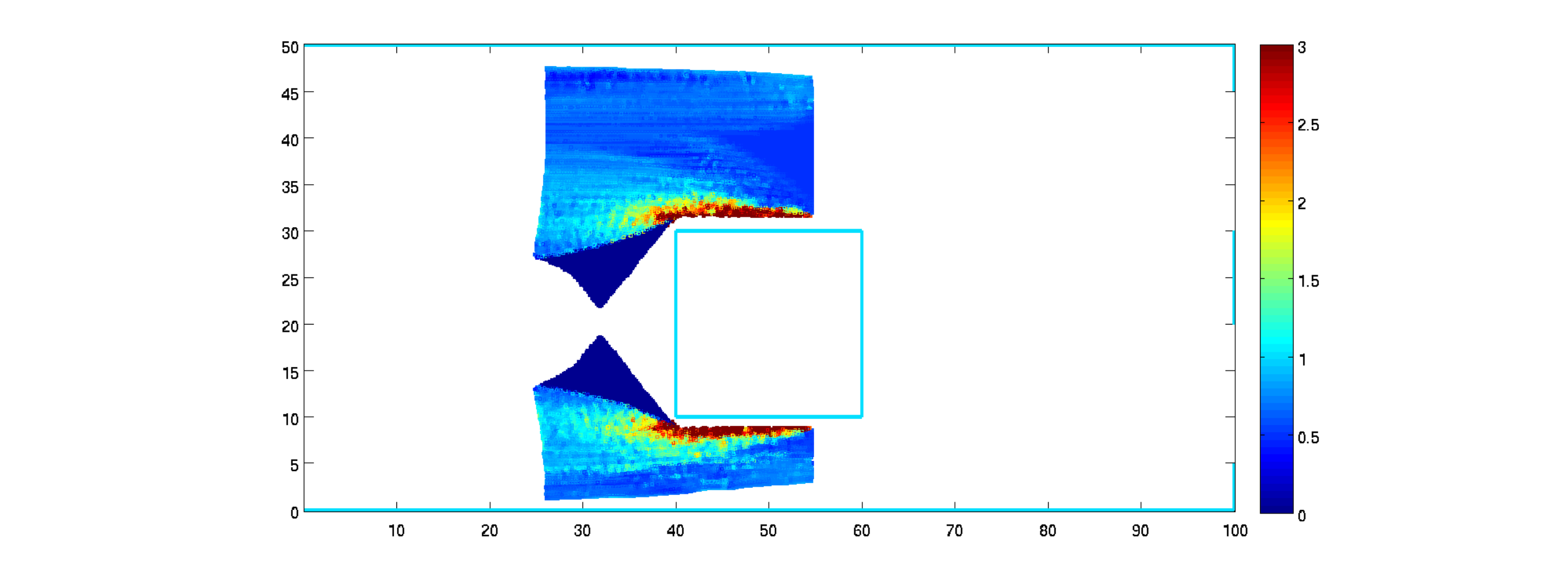}
          	} 
          	\captionsetup[subfigure]{margin=5pt} 
          	\subfloat[$\Delta x =0.2$]{
          		\includegraphics[keepaspectratio=true, width=.5\textwidth]{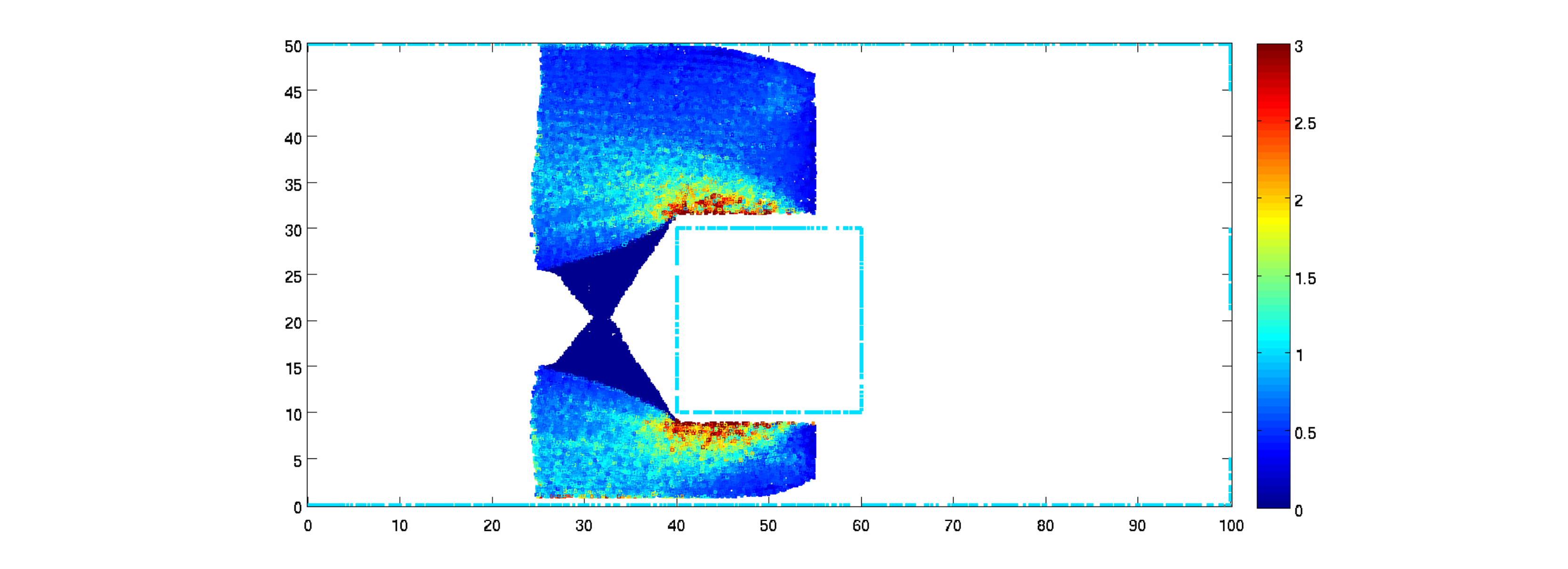}
          	} 
              \caption{Density plot determined from nonlocal equations (\ref{pednonlocalhydro}) with microscopic  approximation at time $t = 12.5$ for $\Delta x = 0.4$ , $\Delta x = 0.3$ and $\Delta x = 0.2$ for $R = 0.2$.   }
             \label{compare_ped1a}
             \end{figure}

    First we consider an underresolved  situation with relatively small 
    value  $R=0.2$ and $\Delta x = 0.5$, i.e a number of particles of approximately $N=5650$.
    In Figure \ref{compare_ped1} we plot the solution at  fixed time using  the localized equation and the nonlocal equation with multiscale method  and  microscopic integration. 
    In this case the solution computed via the multi-scale method 
    and the one computed from  the localized equation coincide, whereas
    the microscopic method gives strongly different results.
   In Figure \ref{compare_ped1a} we show a comparison of solutions obtained from the microscopic scheme with decreasing  discretization sizes or increasing 
   number of particles. 
   
\begin{figure}
   \centering
   \captionsetup[subfigure]{margin=5pt} 
   \subfloat[microscopic]{
   	\includegraphics[keepaspectratio=true, width=.5\textwidth]{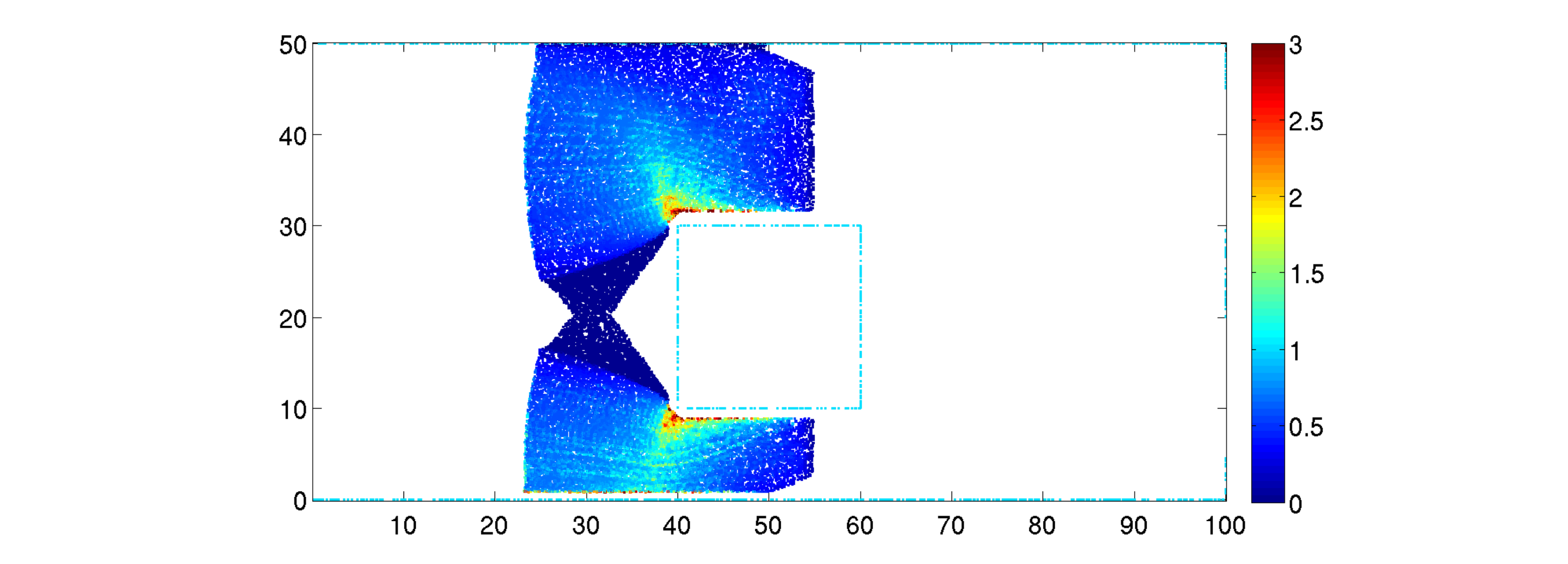}
   } \subfloat[ multi-scale]{
   	\includegraphics[keepaspectratio=true, width=.5\textwidth]{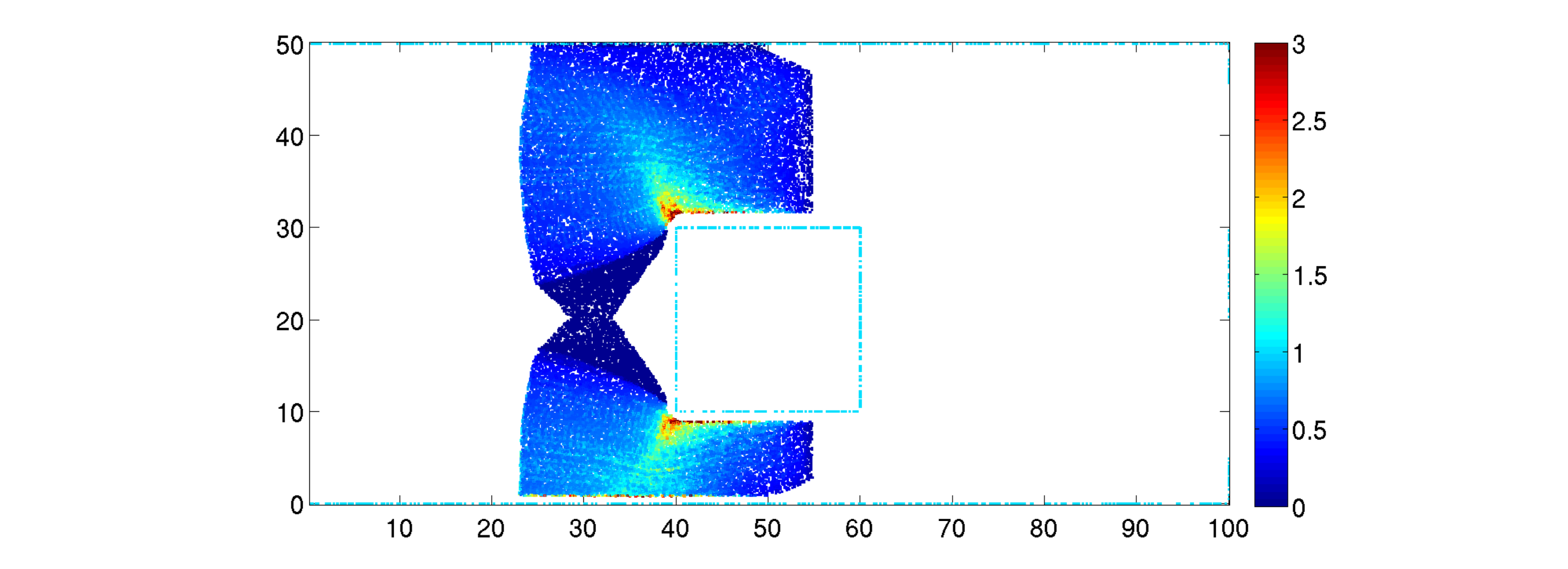}
   } 
      \caption{ Density plot determined from local limit equation (\ref{pedlocalhydro})  and nonlocal equations (\ref{pednonlocalhydro}) with microscopic and multi-scale approximation 
      	 at $t=12,5$ for $\Delta x = 0.2$ and $R=0.4$.}
      \label{compare_ped2}
      \end{figure}

   Figure \ref{compare_ped2} considers a well-resolved case with $\Delta x = 0.2$
   and $R=0.4$. In this case we observe good coincidence of the solutions
   of microscopic an multi-scale approximations.

   Finally, in Table \ref{table3} and Figure \ref{fig115}
   we compute the error and the CPU times of naive and multi-scale  method for different numbers of grid-particles. The errors are determined along a line with $y=37$ and $x \in [25,55]$. The relative $\mathcal{L}^2$-errors  are given  as well as the computation times in minutes.
    The reference solution is computed by  using a spacing of $\Delta x = 0.15$ and
        approximately $62000$ particles. For this  fine resolved case the difference of naive and multi-scale solution is 
        approximately equal to $10^{-2}$.
          Looking at Table \ref{table3} and the multiscale error with $1400$ particles and the naive error with $35200$ particles one observes in this more complex situation    a gain in computation time by more than
     an order of magnitude. comparing  the naive computation with $35200$ particles with the multi-scale simulation with $5700$ 
     particles there is still a gain of  an order of magnitude.


   \begin{figure}
   \centering
      \includegraphics[keepaspectratio=true, angle = 0,   width= .55\textwidth]{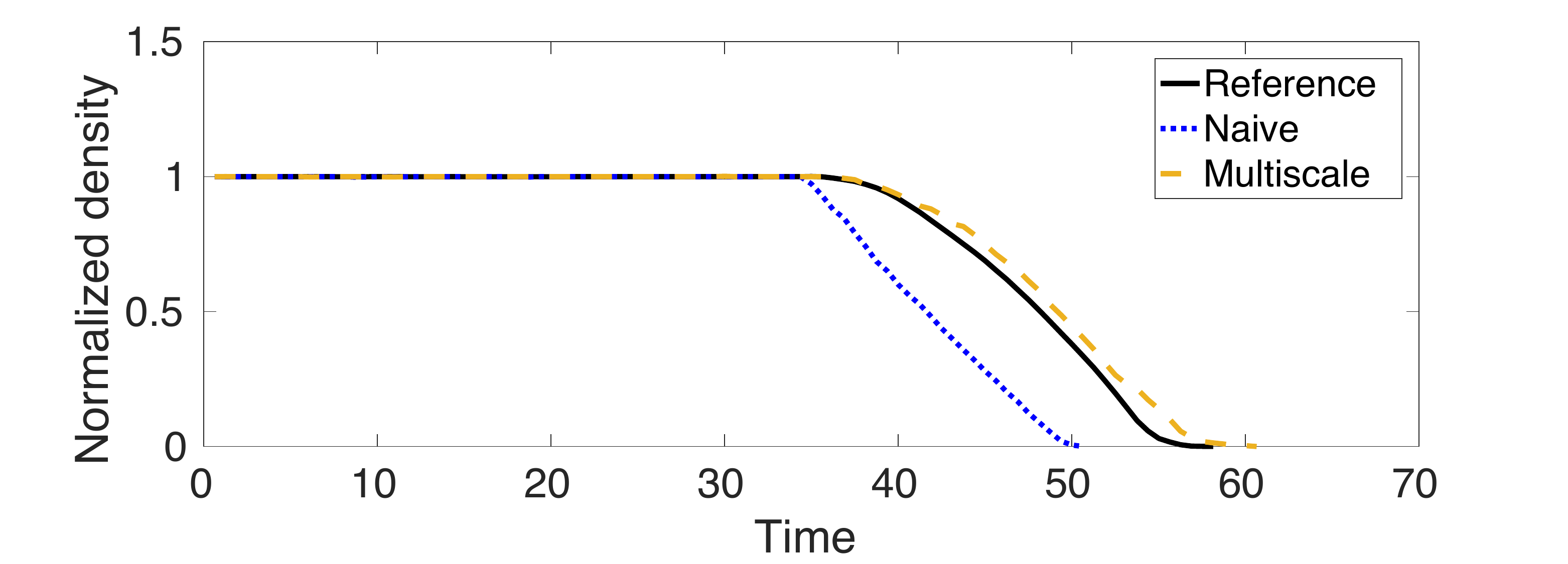}  
      \caption{Time development of the normalized total mass in the computational domain determined from the hydrodynamic pedestrian model  (\ref{pednonlocalhydro}) using  microscopic  and  multi-scale  approximation for fixed interaction radius $R= 0.2$  and coarse  initial spacing $\Delta x =1$ with $N=1400$ grid particles. The reference solution is shown for comparison.}
      \label{fig10}
      \end{figure}

 \begin{table}
 \begin{center}
 \begin{tabular}{|r|r|r|r|r|}
  \hline
 initial   &$\# $ particles  & naive &multi-scale & CPU time  \\ 
  spacing &$ $ &error & error  &  \\ \hline 
 $ 1$ & $1400$ &$ 0.54$ & $0.14 $  & $8$ min   \\
   $0.5$ & $5700$ &$ 0.36$ & $0.18 $  & $23$ min   \\
 $0.35$ &$11500$ &$ 0.48$ &$0.22 $   & $52$ min      \\ 
  $0.2$ & $35200$ & $0.16 $& $0.14 $  & $223$ min  \\   \hline
 \end{tabular}
 \caption{Comparison of CPU times between microscopic and multiscale simulations of the hydrodynamic 2D equations for Example 5. The error analysis is performed at time 12.5 sec. }
 \label{table3}
 \end{center}
 \end{table}

  \begin{figure}
     \centering
        \includegraphics[keepaspectratio=true,   width= .65\textwidth]{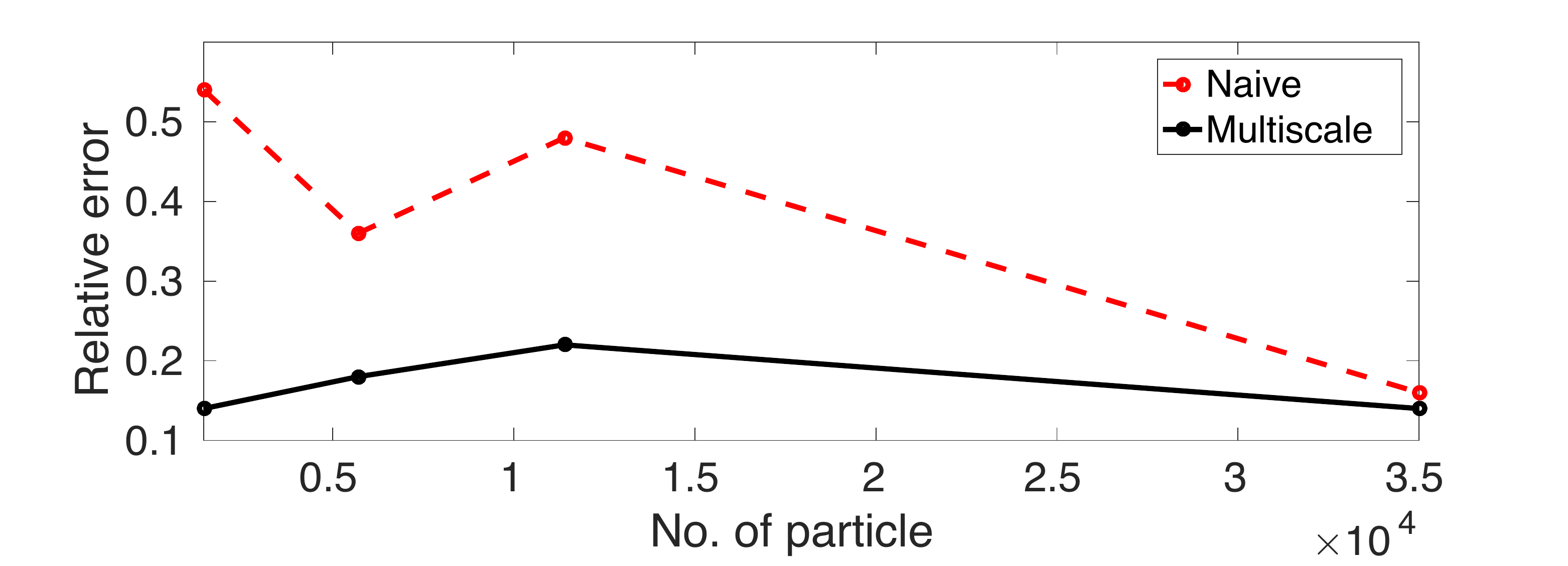}  
        \caption{Error plot for microscopic and multiscale simulations for Example 5.}
              \label{fig115}
        \end{figure}

\section{Conclusion and Outlook }

We have extended a multi-scale meshfree particle method for macroscopic mean field approximations  working uniformly  for  a large range of interaction parameters $R$.
The well resolved case for large $R$ is treated as well as underresolved  situations for small values of $R$.
The method can be considered as a numerical transition from a microscopic system simulation to a macroscopic averaged simulation.
Application of the method are shown for  pedestrian flow simulations. 
The potential  gain in computation time compared to a microscopic simulation
is depending on the situation under consideration wit a potential  gain of several orders of magnitude.
In  the situations considered in this paper we have obtained a gain in computation time of approximately one order of magnitude.
Situations with more complex geometries or  moving boundaries as well as the investigation of attractive-repulsive potentials will be treated in future work, compare \cite{MKT17}.

\subsection*{Acknowledgment}  This  work is supported 
by the German research foundation, DFG grant KL
1105/20-1.



\begin{thebibliography}{10}
\bibliographystyle{siam}


\bibitem{Aro97}
D.G. Aronson, \emph{Regularity properties of flows through porous media}. SIAM J. Appl. Math. 17, (1969), 461-467.

\bibitem{BG15}
S. Blandin and P. Goatin, {\em Well-posedness of a conservation law with non-local flux arising
in traffic flow modeling}. Numer. Math., 132, 2 (2016), 217-241

\bibitem{BH} W. Braun, K. Hepp, \emph{The Vlasov Dynamics
and Its Fluctuations in the 1/N Limit of Interacting Classical
Particles}. Commun. Math. Phys. 56, (1977),  101-113.

\bibitem{BGLX00} 
 T. Belytschko,  Y. Guo, W. Liu, S.P.  Xiao, 
{\em A unified stability analysis of meshless particle methods}
Int. J. Numer. Meth. Engng 48, (2000), 1359-1400.

\bibitem{BLM00}
T. Belyschko, W.K. Liu,B. Moran, 
{\em Nonlinear Finite Elements for Continua and Structures}, John Wiley and Sons, New York, 2000.

\bibitem{BV05}
M. Bodnar, J.L.L. Velazquez, 
\emph{Derivation of macroscopic equations for individual cell-based models: a formal approach},
Math. Meth. Appl. Sci, 28, (2005), 1757-1779 


\bibitem{CP83}
P. Calderoni, M. Pulvirenti,
Propagation of chaos for Burgers' equation
Annales de l'institut Henri Poincaré (A) Physique théorique  39,  1, (1983), 85-97



\bibitem{CCR}
J.A. Ca\~nizo,  J.A. Carrillo, J. Rosado,  \emph{A well-posedness
theory in measures for some kinetic models of collective motion}.
Mathematical Models and Methods in Applied Sciences 21, 3, (2011), 515-539.

\bibitem{CT}
J.A. Carrillo, G. Toscani,  \emph{Asymptotic L1-decay of solutions of the Porous Medium
Equation to self-similarity}. Indiana U. Math. J. 49, (2000), 113-142.



\bibitem{CDP}
J.A. Carrillo, M.R. D'Orsogna, V. Panferov,  \emph{Double milling in
self-propelled swarms from kinetic theory}. Kinetic and Related
Models 2 (2009), 363-378.

\bibitem{CFG}
J.A.  Carillo, M. Di Francesco, M. P. Gualdani,
\emph{Semidiscretization and longtime asymptotics of nonlinear diffusion equations},
Commun. Math. Sci., Supplemental Issue 1, (2007), 21–53.



\bibitem{CFRT}
J.A. Carrillo, M. Fornasier, J. Rosado, G. Toscani, \emph{Asymptotic
Flocking Dynamics for the kinetic Cucker-Smale model},
SIAM J. Math. Anal., 42(1), (2010) 218–236


\bibitem{CGL12}
R. M. Colombo, M. Garavello, and M. Lécureux-Mercier. “A class of nonlocal models for pedestrian traffic”. In: Mathematical Models and Methods in Applied Sciences 22.04 (2012), p. 1150023.

\bibitem{CL12}
R. M. Colombo and M. Lécureux-Mercier. “Nonlocal crowd dynamics models for several populations”. In: Acta Math. Sci. Ser. B Engl. Ed. 32.1 (2012),  177–196.

\bibitem{DM2} 
P. Degond,  F:J: Mustieles,   \emph{Approximation of diffusion equations by deterministic convections of particles}. SIAM J. on Scientific and Statistical Computing 11, (1990),  293-310.


\bibitem{DM1} P. Degond, S. Motsch,
\emph{Continuum limit of self-driven particles with orientation
interaction}. Math. Models Methods Appl. Sci.18  (2008),
1193-1215.

\bibitem{Deg13}
 P. Degond, C. Appert-Rolland, M. Moussaid,
 J. Pettre, G. Theraulaz, 
\emph{A hierarchy of heuristic-based models of crowd
dynamics}, J. Stat. Phys.,  152,  6, (2013), 1033–1068

\bibitem{di}
M. Di Francesco, P.A. Markowich, J.F. Pietschmann and M.T. Wolfram,
\textit{On the Hughes model for pedestrian flow: The one-dimensional case},
J. Differential Equations 250 (2011) 1334-1362.

\bibitem{di2}
M. Di Francesco, S. Fagioli, M.D. Rosini, G. Russo,
Deterministic particle approximation of the {H}ughes model in one space dimension,
Kinetic and Related Models, 10, 1, (2017), 215–237


\bibitem{Dil99}
G. A. Dilts,
{\em Moving-least-squares-particle hydrodynamics: I. Consistency and stability}
 International Journal for Numerical Methods in Engineering 44, (1999), 1115-1155.

\bibitem{DTKB08}
C. Drumm, S. Tiwari, J. Kuhnert, H.-J. Bart,
{\em Finite pointset method for simulation of the liquid-liquid flow field in an extractor}, Comput. chem. Eng. 32, (2008),  2946.

\bibitem{EGKT14}
R. Etikyala, S. Goettlich, A. Klar, S. Tiwari, {\em Particle methods for pedestrian flow models: from microscopic to non-local continuum models}, Mathematical Methods and Models in Applied Sciences 24, 12,
(2014),  2503-2523

\bibitem{FS09}
Y. Farjoun, B. Seibold,
An exactly conservative particle method for one dimensional scalar conservation laws
Journal of Computational Physics, 228, 14, (2009), 5298-5315 




\bibitem{GM77}
R. A. Gingold, J. J. Monaghan,
 \emph{Smoothed Particle Hydrodynamics: theory and application to non-spherical stars}, Mon. Not. Roy. Astron. Soc. 181, (1977), 375-389. 


\bibitem{GS15}
P. Goatin, S. Scialanga,
The Lighthill-Whitham-Richards traffic flow model with non-local velocity: analytical study and numerical results, Netw. Heterog. Media, 11, 1,  (2016), 107-121




\bibitem{MR87} 
S. Mas-Gallic, P. Raviart,
{\em A particle method for first-order symmetric systems} Numerische Mathematik 51, (1987),
323-352.

 
\bibitem{Golse}
F. Golse,
\emph{ On the Dynamics of Large Particle Systems in the mean-field limit}
http://arxiv.org/abs/1301.5494.

\bibitem{HAC74}
C.W Hirt, A.A Amsden, J.L Cook,
{\em An arbitrary Lagrangian-Eulerian computing method for all flow speeds},
Journal of Computational Physics 135, (1997)  203-216


\bibitem{Hel95}
Helbing,D. and P. Molnar, \emph{Social force model for pedestrian dynamics}, Phys. Rev. E, 51
(1995), 4282-4286.


\bibitem{HS09}
J. Haskovec and C. Schmeiser
 Stochastic Particle Approximation for Measure
Valued Solutions of the 2D Keller-Segel System,
 J. Stat. Phys. 135 (2009), 133-151

\bibitem{HS11}
J. Haskovec and C. Schmeiser: Convergence of a stochastic particle approximation
for measure solutions of the 2D Keller-Segel system. Comm. PDE 36 (2011), 940-960.



\bibitem{HPW}
F. Huang, R. Pan, Z. Wang
\emph{L1 Convergence to the Barenblatt Solution
for Compressible Euler Equations
with Damping},
Arch. Rational Mech. Anal. 200 (2011) 665–689.

\bibitem{HS12}
C. Huang, T. W.H. Sheu,
{\em Development of an upwinding particle interaction kernel for simulating incompressible Navier-Stokes equations},Numerical Methods for Partial Differential Equations 28, 5, (2012), 1574–1597.




\bibitem{HMP}
F. Huang, P. Marcati, R. Pan,
\emph{Convergence to the Barenblatt Solution
for the Compressible Euler Equations
with Damping and Vacuum},
Arch. Rational Mech. Anal. 176 (2005) 1–24.

\bibitem{hughes1}
R.L. Hughes,
\textit{The flow of human crowds},
Ann. Rev. Fluid Mech. 35 (2003) 169-182.
 



\bibitem{JX95}
S. Jin and Z. Xin,
 {\em The Relaxation Schemes for Systems of Conservation Laws in Arbitrary Space Dimensions},
Comm. Pure Appl. Math, 48, (1995), 235-277.
   
  
 \bibitem{Jou06} 
 B. Jourdain,
 Probabilistic Approximation via Spatial
Derivation of Some Nonlinear Parabolic
Evolution Equations in 
Monte Carlo and Quasi-Monte Carlo Methods 2004
pp 197-216, Eds. H. Niederreiter, D. Talay, Springer 2006
  
 \bibitem{Jou00} 
 B. Jourdain,
 Diffusion Processes Associated with Nonlinear
Evolution Equations for Signed Measures,
  Methodology and Computing in Applied Probability, 2, 1, (2000),  69-91
 
 
 \bibitem{KT15}
 A. Klar, S. Tiwari,  A multi-scale  meshfree particle method for macroscopic mean field  interacting particle models, SIAM Multiscale Mod. Sim. 12, 3, (2014), 1167–1192
 
 
\bibitem{LWZSL09} 
 H. Ling, S. C. Wong, M. Zhang, C. H. Shu and W. H. K. Lam, Revisiting Hughes
dynamic continuum model for pedestrian flow and the development of an efficient
solution algorithm, Transp. Res. B, Methodol. 43, (2009), 127–141
 
 \bibitem{MKT17}
 N. Mahato, A. Klar, S. Tiwari, Modelling and simulations of macroscopic multi-group pedestrian flow,  in preparation
 
 
  \bibitem{SK06} 
  A. Sopasakis and M. A. Katsoulakis. “Stochastic modeling and simulation of traffic flow: asymmetric single exclusion process with Arrhenius look-ahead dynamics”. In: SIAM J. Appl. Math. 66.3 (2006), 921–944 (electronic).
  
  
  \bibitem{LTB09} 
  A. J. Leverentz, C. M. Topaz, and A. J. Bernoff,
  {\em Asymptotic Dynamics of Attractive-Repulsive Swarms},
  SIAM J. Appl. Dyn. Systems 8,  3, (2009), 880-908

  
    
\bibitem{neunzert}
H. Neunzert, \emph{The Vlasov equation as a limit of Hamiltonian
classical mechanical systems of interacting particles}. Trans.
Fluid Dynamics, 18 (1977), 663-678.

\bibitem{oehl}
K. Oelschlaeger, \emph{Large systems of interacting particles and the porous medium equation}. Journal of Differential
Equation 88(2), (1990), 294–346

\bibitem{scho}
M.E. Schonbeck, \emph{Convergence of solutions to nonlinear dispersive equations}, Comm. Partial Differ. Eq., 7, (1982) , 959-1000 


\bibitem{spohn1}
H. Spohn, \emph{Kinetic equations from Hamiltonian dynamics: Markovian limits}, Rev. Mod. Phys. 52,  (1980) 569–615 


\bibitem{spohn2}
H. Spohn,  \emph{Large scale dynamics of interacting particles}. Texts
and Monographs in Physics, Springer (1991).

\bibitem{szn}
A.S. Sznitman,\emph{A propagation of chaos result for Burgers' equation}
Prob. Th. Rel. Fields 71, (1986), 581-613 
+
\bibitem{TK07}
 S. Tiwari, J. Kuhnert, \emph{Modelling of two-phase flow with surface tension by Finite Point-set method (FPM)}.
 J. Comp. Appl. Math. 203 (2007),  376-386.

\bibitem{TK03}
S. Tiwari, J.  Kuhnert, 
{\em Finite pointset method based on the projection method for simulations of the incompressible Navier-Stokes equations},
 M. Griebel, M. A. Schweitzer (Eds.), Springer LNCSE: Meshfree Methods for Partial Differential Equations, Springer-Verlag, Berlin 26, (2003),  373-387.
 
\bibitem{TKH09}
S. Tiwari, A. Klar, S. Hardt, 
\newblock \emph{A particle-particle hybrid method for kinetic and 
continuum equations}, J. Comp. Phys. 228, (2009), 7109-7124. 


\bibitem{TTW162}
D. R. Tunuguntla,  A.R. Thornton, T. Weinhardt, 
\emph{From discrete particles to continuum fields: extension to bidisperse systems}
Comput. Part. Mech. 3,3, (2016), 349-365

\bibitem{WTLB12}
T. Weinhardt, A.R. Thornton, S. Luding, O. Bokhove, 
\emph{From discrete particles to continuum fields near a boundary}
Granul. Matt. 14,2, (2012), 289-294

\bibitem{ZY02}
H.P. Zhu, A.B. Yu, \emph{Averaging method of granular materials}. Phys. Rev. E 66 (2002), 021302.



  \end{thebibliography}
\end{document}